\documentclass{article}
\usepackage[top=1in, bottom=1in, left=1in, right=1in]{geometry}
\usepackage{amsmath,amsthm,amssymb,amsfonts,mathtools,array,longtable,float,booktabs}
\usepackage{caption}
\usepackage{graphicx,subcaption}
\usepackage[colorlinks=true, pdfstartview=FitV,linkcolor=blue,citecolor=blue,urlcolor=blue]{hyperref}
\usepackage{appendix}
\usepackage{todonotes}

\usepackage{mathrsfs,slashed}

\let\EellsSin\sin

\let\EellsMathOpen\(
\let\EellsMathClose\)

\newcommand{\Z}{\mathbb{Z}}

\newcommand{\R}{\mathbb{R}}
\newcommand{\C}{\mathbb{C}}

\DeclareMathOperator{\ind}{index}

\renewcommand{\Im}{\operatorname{Im}}
\renewcommand{\Re}{\operatorname{Re}}

\renewcommand{\epsilon}{\varepsilon}
\DeclareMathOperator{\sen}{sen}
\renewcommand{\sin}{\sen}

\def\({\mathopen{}\left(}
\def\){\right)\mathclose{}}
\def\<{\mathopen{}\left<}
\def\>{\right>\mathclose{}}

\usepackage{multicol,color}

\definecolor{gold}{rgb}{0.85,.66,0}
\definecolor{cherry}{rgb}{0.9,.1,.2}
\definecolor{burgundy}{rgb}{0.8,.2,.2}
\definecolor{orangered}{rgb}{0.85,.3,0}
\definecolor{orange}{rgb}{0.85,.4,0}
\definecolor{olive}{rgb}{.45,.4,0}
\definecolor{lime}{rgb}{.6,.9,0}
\definecolor{green}{rgb}{.2,.7,0}
\definecolor{grey}{rgb}{.4,.4,.2}
\definecolor{brown}{rgb}{.4,.3,.1}

\usepackage{xparse}

\makeatletter
\NewDocumentCommand{\NewExampleTag}{mm}{%
  \expandafter\newcommand\csname tag#1\endcsname[1]{%
    \ifdefined\mmriemexercisereplayactive\else
      \addcontentsline{explist#1}{#1}{##1}%
    \fi
    \ifdefined\mmriemexerciseactive
      \textnormal{(cf.~\ref{exemplos:#2}).}\space
    \fi
  }%
  \expandafter\newcommand\csname showtag#1\endcsname{%
    \section{#2}\label{exemplos:#2}%
    \@starttoc{explist#1}%
  }%
  \expandafter\newcommand\csname l@#1\endcsname[2]{%
    \textbf{##1} -- pág.~##2\par
  }%
}
\makeatother

\NewExampleTag{euclid}{Espaço euclidiano}
\NewExampleTag{esfera}{Esfera}
\NewExampleTag{reproj}{Espaço projetivo real}
\NewExampleTag{comproj}{Espaço projetivo complexo}
\NewExampleTag{regval}{Variedades de nível}
\NewExampleTag{superficie}{Superfícies regulares}
\NewExampleTag{homog}{Espaços homogêneos}
\NewExampleTag{lie}{Grupos de Lie}

\let\sin\EellsSin
\let\(\EellsMathOpen
\let\)\EellsMathClose

\let\C\relax
\usepackage{constants}
\let\EellsConstant\C
\let\C\EellsComplex
\renewcommand{\Cl}[2][normal]{\EellsConstant[#1]\labelconstant{#2}{\string #1}}
\usepackage[nameinlink,noabbrev]{cleveref}
\newcommand{\ghlogo}{\raisebox{-0.2em}{\includegraphics[height=1em]{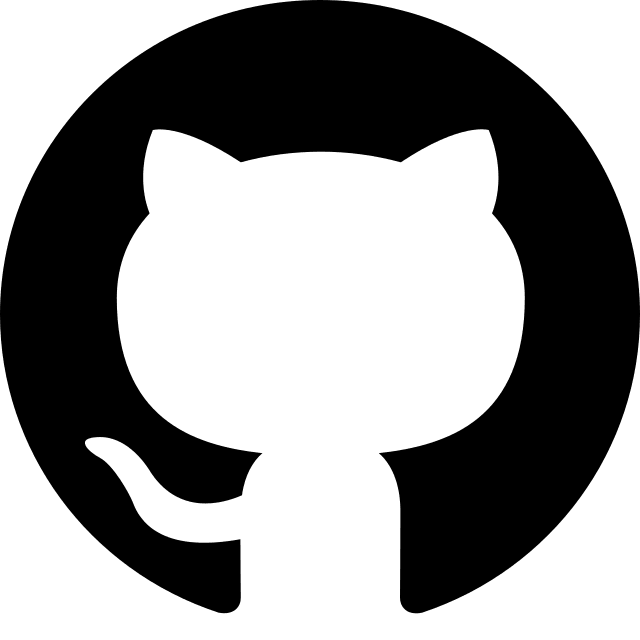}}}

\newtheorem{theorem}{Theorem}[section]
\newtheorem{conjecture}{Conjecture}[section]
\crefname{conjecture}{Conjecture}{Conjectures}

\renewcommand{\thefootnote}{\fnsymbol{footnote}}

\newtheorem{lemma}[theorem]{Lemma}
\newtheorem{proposition}[theorem]{Proposition}

\theoremstyle{definition}
\newtheorem{remark}[theorem]{Remark}
\newtheorem{definition}[theorem]{Definition}
\theoremstyle{plain}

\title{Numerical verification for harmonic maps:\\ Counterexamples to Eells' conjecture}
\author{Edward Hirst\textsuperscript{1}\thanks{\href{mailto:ehirst@unicamp.br}{\texttt{ehirst@unicamp.br}}},
Daniel Platt\textsuperscript{2}\thanks{\href{mailto:d.platt@imperial.ac.uk}{\texttt{d.platt@imperial.ac.uk}}},\\
Henrique N. S\'a Earp\textsuperscript{1}\thanks{\href{mailto:hqsaearp@unicamp.br}{\texttt{hqsaearp@unicamp.br}}},
Jakob Stein\textsuperscript{1}\thanks{\href{mailto:jstein@unicamp.br}{\texttt{jstein@unicamp.br}}}}
\date{\small\today \\[0.5em]
\begin{minipage}{0.94\textwidth}\centering\small\itshape
\textsuperscript{1}Instituto de Matem\'atica, Estat\'istica e Computa\c{c}\~ao Cient\'ifica (IMECC), \\
Universidade Estadual de Campinas (UNICAMP), 13083-859, Brasil\\[0.5em]
\textsuperscript{2}Department of Mathematics, Imperial College London, \\ 
180 Queen's Gate, South Kensington, London SW7 2RH, United Kingdom
\end{minipage}}

\begin{document}
\maketitle
\renewcommand{\thefootnote}{\arabic{footnote}}
\setcounter{footnote}{0}
\numberwithin{equation}{section}
\numberwithin{figure}{section}
\numberwithin{table}{section}
\setcounter{tocdepth}{3}

\begin{abstract}
This paper develops a framework to construct harmonic maps $S^3 \rightarrow S^2$ numerically, which is used to construct approximate harmonic maps in homotopy classes with Hopf invariant $2,3,4,5,6$. By using a rigorous fixed-point argument, it is shown that there exist genuine harmonic maps with the same Hopf invariant close to these approximate solutions, thus producing a family of counterexamples to Eells' conjecture, and the related Ding's conjecture.\\
\small Code: \ghlogo \, \href{https://github.com/edhirst/eells_counterexamples}{\texttt{github.com/edhirst/eells\_counterexamples}}
\end{abstract}

\clearpage
\tableofcontents
\newpage

\section{Introduction}
\label{sec:intro}

Harmonic maps $u: M \rightarrow N$ are defined as critical points of the Dirichlet energy functional
\[
E(u)=\frac12\int_M |du|^2\,d\mathrm{vol}_M.
\]
There are well-known obstacles to finding such critical points:
first, often, the infimum of the energy in a fixed topological class is zero (see \cite[Section~2.5]{Struwe} for a precise statement).
Hence, direct minimisation methods for $E$ cannot produce harmonic maps in these cases.
Second, in many interesting situations, harmonic maps are proven to not even be a \emph{local} minimum, but their second variation has a positive index (see e.g. \cite[Theorem~I.1]{Riviere}).

The result is a large body of research constructing examples of harmonic maps and analysing their obstructions using sophisticated methods:
the heat flow method from \cite{EellsSampson}, the min-max constructions in \cite{SacksUhlenbeck,KarpukhinStern2024a,KarpukhinStern2024b}, using symmetries and ODE methods in \cite{Smith,EellsRatto,Bizon}, the gluing constructions in \cite{BrezisCoron,Jost} and gluing obstructions in \cite{Rupflin}, and the holomorphic and twistor constructions in \cite{EellsWood,EellsSalamon}, to list a few.

On the other hand, a great variety of numerical techniques has been applied to numerically compute harmonic maps.
An excellent overview over different approaches is provided in \cite{BBPS}.
The main approaches include:
finite element methods based on the energy-decreasing iteration
of \cite{Alouges}, discretised in \cite{Bartels}, and the more
recent approximation results for stable harmonic maps in
\cite{BPW}, or approximating the harmonic-map heat flow in
\cite{BKW};
finite-difference methods, e.g. the recent \cite{XWXZ};
lattice discretisations \cite{MullerStruweSverak} and intrinsic geometric methods such as \cite{GasterLoustauMonsaingeonDiscrete,GasterLoustauMonsaingeon};
the approach using point clouds from \cite{WuYau};
and ODE techniques for symmetry-reduced solutions in \cite{BizonChmaj}.

An important special case that has received much attention is the case of harmonic maps between round spheres, and among these the case of $u:S^3 \rightarrow S^2$ in particular.
In this situation, nonconstant solutions must be unstable by \cite[Theorem~I.1]{Riviere}.
The Hopf fibration $\pi:S^3 \rightarrow S^2$ is a classical example of a harmonic map in this case.
More generally, one checks that $R \circ \pi \circ A$ is a harmonic map for $R: S^2 \rightarrow S^2$ a weakly conformal map and $A: S^3 \rightarrow S^3$ an isometry \cite{BairdWood}.
To our knowledge, no other harmonic maps $S^3 \rightarrow S^2$ have appeared in the literature, neither with a mathematical proof, nor using numerical experiments.
Indeed, Eells' conjecture (see \cite[p.~730]{Wang}) says:

\begin{conjecture}[Eells]
\label{conj:eells}
Every harmonic map $u:S^3 \rightarrow S^2$ has the form
$u=R \circ \pi \circ A$ for a weakly conformal map
$R:S^2 \rightarrow S^2$ and an isometry $A:S^3 \rightarrow S^3$.
\end{conjecture}

To a smooth map $u: S^3 \rightarrow S^2$ one can assign its (absolute) \emph{Hopf invariant}, see \S \ref{subsection:hopf-invariant}.
The Hopf invariant of the Hopf fibration $\pi$ is $1$, and more generally for a composition $R \circ \pi \circ A$ it is $|\deg A|(\deg R)^2$ \cite[Section~2]{Wang}.
\Cref{conj:eells} implies the following weaker conjecture,
attributed to Ding in \cite[p.~730]{Wang}:

\begin{conjecture}[Ding]
\label{conj:ding}
A homotopy class of maps $S^3 \rightarrow S^2$ contains a harmonic map
if and only if its (absolute) Hopf invariant is $k^2$ for some integer $k$.
\end{conjecture}

The content of the present paper is the following:

\begin{itemize}
\item
Rewriting the defining equation for harmonic maps $u:S^3 \rightarrow S^2$ as a polynomial equation for maps $u:S^3 \rightarrow \R^3$, we can take an ansatz for $u$ as a linear combination of spherical harmonics, and use spectral methods to solve the equation.
Apart from the usual good convergence properties of spectral methods, we find that the method is robust and for many initial data converges to approximate harmonic maps with any given Hopf invariant.
In particular, we exhibit numerical examples that appear to contradict \cref{conj:eells,conj:ding}, though this does not constitute a proof.

\item
We then apply a fixed-point argument using elliptic PDE techniques to prove that near each considered numerical example, there exists a genuine harmonic map.
We quantify exactly what \emph{near} means, and are able to compute its Hopf invariant and tight enclosures for its energy.
This turns each numerical candidate from the previous bullet point into a proven counterexample to the conjectures above.
\end{itemize}

\newpage
Our main result is:

\begin{theorem}
\label{thm:main}
	There exist harmonic maps $u:S^3 \rightarrow S^2$ with the Hopf invariants and energies listed in \Cref{tab:energies}. None of the maps are of the form $u=R \circ \pi \circ A$ for an isometry $A: S^3 \rightarrow S^3$ and a weakly conformal map $R: S^2 \rightarrow S^2$.
\end{theorem}

\begin{table}[ht] 
\centering
\begin{tabular}{@{}lccccc@{}}
\toprule
$H(u)$ & 2 & 3 & 4 & 5 & 6 \\
\midrule
$E(u)/8\pi^2$ & 1.587 & 1.971 & 2.251 & 2.472 & 2.652 \\
\bottomrule
\end{tabular}
\caption{Hopf invariants and normalised Dirichlet energies $\tfrac{E}{8\pi^2}$ of our harmonic maps $u$, rounded to three decimal places.}
\label{tab:energies}
\end{table}
In particular, this disproves \cref{conj:ding} and thereby the stronger \cref{conj:eells}.
It also disproves \cref{conj:eells} for cases with square Hopf invariant.

The technique of constructing an approximate solution to an elliptic PDE and perturbing it to a genuine solution has been a standard tool in geometric analysis since \cite{Taubes} and has been applied in countless important works.
These proofs usually use pen-and-paper techniques.
Much more recently, this same approach has been applied in geometry using computer aid in the form of \emph{numerically verified proofs}.
In \cite{ButtsworthHodgkinson,WangSU2,SalamonSinghal,CaoLaboraRodriguezVazquez}, this was carried out using symmetries, reducing the analysis to an ODE.
In \cite{Cortes:2026kfx, PlattNirenberg}, this was applied to a PDE on the $2$-sphere.

In our case, we also impose symmetries on our approximate solution and all function spaces, which guarantee that the genuine solution has the same symmetries.
We consider a family of actions of $G=O(2) \times \mathbb{Z}_2$ on $S^3$ and $S^2$ (see \cref{subsection:sobolev-spaces}), searching for equivariant maps.
These actions have cohomogeneity two on the domain, so we cannot apply the ODE techniques from the references above, but apply general PDE techniques that work in any dimension.
However, the reduction to two dimensions provides important speed-ups for our numerical work. 

\begin{remark} The harmonic maps in \Cref{conj:eells} have energy $E = 8 \pi^2 |\deg{R}| = 8 \pi^2 \sqrt{H}$, so the energies listed in \Cref{tab:energies} exclude our solutions in \cref{thm:main} from being of this form. Moreover, at least when the Hopf invariant is even, one can a-priori exclude any map of the form in \cref{conj:eells} from our equivariant function spaces, see Remark \ref{rem:equivariantq}. 
\end{remark}

This proof method requires to give precise bounds for all quantities involved:
residuals, linear operators, and higher order terms.
It requires reducing analysis in infinite-dimensional Banach spaces to finite computations that can be implemented on a computer.
While applications in geometry are still few, we build on a solid foundation of existing numerically verified proofs in analysis.
We mention the important existence proof of the Lorenz attractor in \cite{Tucker} as an example, and refer to the survey \cite{GomezSerrano} and the textbook \cite{Nakao} for comprehensive overviews.

One important feature of the proof technique is that we obtain very detailed information about the genuine solutions coming from the approximate solutions.
Apart from the energies shown in \cref{thm:main}, we demonstrate this principle with the following result:

\begin{theorem}
    \label{thm:index}
    The harmonic map $u:S^3 \rightarrow S^2$ with Hopf invariant $2$ from \cref{thm:main} has $\ind(u) \geq 9$.
\end{theorem}

There are $4$ distinguished vector fields on $S^3$ that are given as the gradients of the four coordinate functions. Precomposing $u$ with their flows gives negative directions for the second variation of the Dirichlet energy, see \cite[Theorem~1.1]{ElSoufi} and \cite[\S II.5]{Riviere}.
Thus, $\ind(u)\ge4$.
Even more: \cite[Theorem~I.1]{Riviere} proves $\ind(u)\ge5$ if $u$ is not of the form $R\circ\pi\circ A$, and our lower bound strengthens the result in this example.

The paper is structured as follows:
in \cref{sec:background} we make the explicit choices that are used in the later analysis.
In particular, in \cref{subsection:sobolev-spaces} we introduce the spherical harmonics and in \cref{subsection:harmonic-map-equation} we rewrite the harmonic map equation as a \emph{polynomial} equation.
Following this, we prove the estimates for the residual, linearised operator, and higher order terms in \cref{section:estimates}.
The final section \cref{section:conclusion} is then a simple matter of assembling the earlier estimates to obtain proofs for \cref{thm:main,thm:index}.
The appendices list the central data used in our certified computations.

\section{Background}
\label{sec:background}

In this section we marshal some definitions and notations that are used throughout the remainder of the paper.
In \cref{subsection:sobolev-spaces}, we define the Sobolev spaces of equivariant maps that we will use for our functional analysis.
In \cref{subsection:harmonic-map-equation}, we rewrite the harmonic map equation for polynomials.
Following this, in \cref{subsection:hopf-invariant}, we define the Hopf invariant and prove a simple proposition that will later allow us to compute the Hopf invariant of our approximately harmonic candidate map using a numerical bound. In \cref{subsection:energy}, we compute the energy of the harmonic maps of \cref{conj:eells}, in order to compare to the energy of our solutions. In \cref{subsection:index}, we introduce the quadratic forms and estimates used to obtain the index lower bound. Last, in \cref{subsection:fixed-point} we will state a fixed-point theorem that will be applied at the very end of the paper to perturb the approximately harmonic map to a genuine harmonic map.

\subsection{Sobolev spaces with symmetries}
\label{subsection:sobolev-spaces}

In this subsection, we define Banach spaces of functions that have the symmetries we expect of our solution.
We will then carry out all functional analysis in these function spaces and show that a solution indeed exists in them.
This approach of defining symmetry-reduced function spaces is standard for elliptic PDEs.
An example for harmonic maps is \cite{CorletteWald}.
In our case, we will consider \emph{equivariant} maps $S^3 \rightarrow S^2$, hence we begin by defining group actions on $S^3$ and $S^2$, or, more precisely, on the ambient $\R^4$ and $\R^3$ that restrict to the spheres.

\subsubsection{Definition of the symmetry group.}
\label{sec:symmetries}
Identify $\R^4=\C^2$ and $\R^3=\C\times\R$.
Let $G=O(2)\times\Z_2$, with circle subgroup $e^{i\theta}$, reflection $\kappa\in O(2)$
and generator $\epsilon\in\Z_2$. 
For $z,w,\zeta \in \C$ and $t \in \R$, and some positive integer $q>1$, we define the action of $G$ on $\R^4$ and $\R^3$ by the following,

\begin{align} \label{eqn:symmetries}
\begin{aligned}
e^{i\theta}\cdot(z,w)&=(e^{i\theta}z,e^{-qi\theta}w),
&e^{i\theta}\cdot(\zeta,t)&=(e^{-2qi\theta}\zeta,t),\\
\kappa\cdot(z,w)&=(\bar z,\bar w),
&\kappa\cdot(\zeta,t)&=(\bar\zeta,t),\\
\epsilon\cdot(z,w)&=(z,-w),
&\epsilon\cdot(\zeta,t)&=(\zeta,-t).    
\end{aligned}
\end{align}
All these actions restrict to the round unit spheres $S^3$ and $S^2$. We also define the following spaces of equivariant functions:
\begin{align*}
X&:=H^2_G(S^3,\R^3):=\{u\in H^2(S^3,\R^3):u(gx)=g u(x)\text{ for all }g,x\},\\
Y&:=L^2_G(S^3,\R^3):=\{u\in L^2(S^3,\R^3):u\circ g=g u\text{ in }L^2\text{ for all }g\in G\},
\end{align*}
and endow them with the norms $\|u\|_X=\|(1-\Delta)u\|_{L^2}$ and $\|u\|_Y=\|u\|_{L^2}$. Here, and throughout the text, we will refer to the group $G$ as $G_q$ when we want to specify that it acts via \eqref{eqn:symmetries} for some given $q$. Otherwise, we will just refer to $X$, $Y$ and $G$ when our analysis is independent of $q$. 
\begin{remark} \label{rem:equivariantq} When $q$ is even, the action of $(e^{i\pi},\varepsilon) \in G_q$ in \eqref{eqn:symmetries} is given by 
\begin{align*}
(z,w) \mapsto (-z,-w)& &(\zeta, t)\mapsto (\zeta,-t) 
\end{align*}
This excludes a $G_q$-equivariant map $u: S^3 \rightarrow S^2$ being of the form $R \circ \pi \circ A$ , where $R:S^2 \mapsto S^2$ is surjective, $\pi:S^3 \rightarrow S^2$ is the Hopf projection, and $A \in O(4)$. One can argue as follows: since $\pi(-z,-w) = \pi  (z,w)$ and $A(-z,-w) = -A(z,w)$ by linearity, then $u(z,w) =R \circ \pi \circ A (z,w) = u (-z,-w)$. If $u$ is also equivariant, then the image of $u$ must lie in the plane $\C$ fixed by $(\zeta, t)\mapsto (\zeta,-t)$, so $R$ cannot be surjective. In particular, for $q$ even, this rules out harmonic maps of this form occurring in our function spaces $X$ and $Y$. Indeed, our previous argument shows that $u$ would be a harmonic map $S^3 \rightarrow S^1$, which must be constant. 
\end{remark}
\begin{remark} The Hopf fibration is equivariant with respect to the action of \eqref{eqn:symmetries} with $q=1$, see \cref{def:hopffibration}, so the restriction to $q >1$ is important for our construction. 
\end{remark}
\subsubsection{Spherical harmonics.}
Throughout the non-rigorous initial numerical work as well as for the later functional analysis, we crucially use the spherical harmonics on $S^3$.
During optimisation, they allow us to use the spectral methods described in \cref{subsection:computing-the-approx-sol}.
For functional analysis, they allow us to decompose the linearised elliptic operator into blocks that can be computed effectively.
We use the non-positive Laplacian $\Delta=\operatorname{div}\nabla$.
The complex spherical harmonics of degree $l$ are the restrictions to $S^3$ of
complex homogeneous harmonic polynomials of degree $l$ on $\R^4$.
They form the $(l+1)^2$-dimensional eigenspace of $-\Delta$ with eigenvalue $l(l+2)$.
We use the Jacobi basis of \cite[\S2, (4)]{LachiezeRey}, normalized for volume $2\pi^2$:
\begin{align}
Y_{m,n,j}(z,w)&=c_{|m|,|n|,j}\,z^{(m)}w^{(n)}P_j^{(|m|,|n|)}(1-2|z|^2),\\
c_{a,b,j}^2&=\frac{(2j+a+b+1)(j+a+b)!\,j!}{2\pi^2(j+a)!(j+b)!},
\qquad l=|m|+|n|+2j.
\end{align}
Here $m,n\in\Z$, $j\ge0$, and
$z^{(m)}=z^m$ for $m\ge0$, $z^{(m)}=\bar z^{|m|}$ for $m<0$, similarly for $w^{(n)}$. Moreover $P_j^{(a,b)}$ is the Jacobi polynomial with the real variable $(1-2|z|^2)$, i.e. 
\begin{align*}
  P_j^{(a,b)}(1-2|z|^2) = \sum_{s=0}^{j}(-1)^s\binom{j+a}{j-s}\binom{j+a+b+s}{s}\,|z|^{2s}. 
\end{align*}
For each degree $l$, the functions $Y_{m,n,j}$ with $(m,n,j)$ such that  $|m|+|n|+2j=l$ form an $L^2$-orthonormal basis of the Laplacian eigenspace. 
\begin{lemma} \label{lem:action}The action of $G_q$ on $Y_{m,n,j}(z,w)$ is as follows:
\begin{enumerate}
    \item $e^{i\theta} \cdot Y_{m,n,j} = e^{i(m - qn)\theta} Y_{m,n,j}$
    \item $\kappa \cdot Y_{m,n,j} = Y_{-m,-n,j} = \overline{Y_{m,n,j}}$
    \item $\epsilon \cdot Y_{m,n,j} = (-1)^{n}Y_{m,n,j}$
\end{enumerate}
\end{lemma}
We now specialise to the \textit{equivariant} subset of our basis, where the target $\R^3$ is now also equipped with a $G$-action \eqref{eqn:symmetries}. 

\begin{proposition}\label{prop:harmonic-bases}
Let $G =G_q$. The following fields are $L^2$-orthonormal and their linear spans are dense in $X$ and $Y$ with respect to their respective norms:
\begin{align*}
&(\Re Y_{q(n-2),n,j},\Im Y_{q(n-2),n,j},0),&&n \; \text{even},\quad j\ge0,\\
&(0,0,\sqrt2 \Re Y_{qn,n,j}),&&n\ge1\text{ odd},\quad j\ge0.
\end{align*}
\end{proposition}

\begin{proof}
We identify $\mathbb{R}^3 = \mathbb{C} \times \mathbb{R}$ as before. The $G$-action \eqref{eqn:symmetries} preserves the components in this splitting, and we determine which members of the Jacobi basis $Y_{m,n,j}$ are equivariant for each component. In the first component $\mathbb{C}$, part (1) of Lemma \ref{lem:action} gives $m-qn = -2q$, and part (3) gives that $n$ must be even. The second component $\mathbb{R}$ is similar, but now only $\epsilon \in G$ acts non-trivially on the target (as $-1$), so we must have $m=qn$, and $n$ odd. The restriction $n\geq 1$ avoids the overcounting coming from (2), since $\Re Y_{qn,n,j} = \Re Y_{-qn,-n,j}$.
\end{proof}

Enumerate the vector basis by $(\Psi_k)$ and denote its degrees by $l_k$. Then for some $(m_k,n_k,j_k)$ we have $l_k = |m_k| + |n_k|+2j_k$, which gives $l_k \geq 2$ by Proposition \ref{prop:harmonic-bases}. Moreover, we have that
\begin{align}
\label{eqn:X-norm-calculation}
\Big\|\sum_k a_k\Psi_k\Big\|_X^2=\sum_k(l_k+1)^4|a_k|^2.
\end{align}

\subsubsection{Sobolev norms.}
In addition to the $X$-norm, we use the Sobolev norm
$\|\xi\|_{H^1}^2
=\|\xi\|_{L^2}^2+\|\nabla\xi\|_{L^2}^2$.
We have the following Sobolev embedding theorems for the chosen Sobolev norms:

\begin{proposition}\label{prop:sobolev}
For $u\in X$ and $\xi\in H^1(S^3,\mathbb R^3)$,
\[
\|u\|_{C^0}\le c_\infty\|u\|_X,\quad
\|\nabla u\|_{L^4}\le c_4\|u\|_X,\quad
\|\Delta u\|_{L^2}\le\|u\|_X,
\quad
\|\xi\|_{L^4}\le\beta\|\xi\|_{H^1},
\quad
\|\xi\|_{L^4}\|\nabla\xi\|_{L^2}\le\beta\|\xi\|_{H^1}^2
\]
where
\[
c_\infty=\sqrt{\frac{\pi^2/6-5/4}{2\pi^2}},
\qquad
c_4=\left(\frac{2}{9\pi^2}\right)^{1/4},
\qquad
\beta=(2\pi^2)^{-1/4}.
\]
\end{proposition}

\begin{proof}
It suffices to consider finite harmonic sums
$u=\sum_{l\ge2}u_l$, where $u_l$ has degree $l$.
The addition formula \cite[Corollary~1.2.7]{DaiXu},
applied componentwise, and Cauchy--Schwarz give
\[
\begin{aligned}
\|u\|_{C^0}
&\le \frac1{\sqrt{2\pi^2}}
\sum_{l\ge2}(l+1)\|u_l\|_{L^2}
\le
\left(\frac1{2\pi^2}\sum_{l\ge2}\frac1{(l+1)^2}\right)^{1/2}
\|u\|_X
=c_\infty\|u\|_X.
\end{aligned}
\]
The Laplacian estimate follows from $l(l+2)\le(l+1)^2$.

The sharp interpolation inequality on $S^3$ \cite[Eqn. 6.8.7]{BakryGentilLedoux},
after rescaling the measure, gives
\[
\|f\|_{L^4}^2
\le \frac1{\sqrt{2\pi^2}}
\left(\|f\|_{L^2}^2+\frac23\|\nabla f\|_{L^2}^2\right)
\qquad(f\in H^1(S^3)),
\]
which immediately gives the two claims for $\xi$.
Applying this to $f=|\nabla u|$, Kato's inequality and the
integrated Bochner formula on $S^3$ give
\[
\begin{aligned}
\|\nabla u\|_{L^4}^2
&\le
\frac{\|\nabla u\|_{L^2}^2
+\frac23\|\nabla^2u\|_{L^2}^2}{\sqrt{2\pi^2}}
=
\frac{2\|\Delta u\|_{L^2}^2-\|\nabla u\|_{L^2}^2}
{3\sqrt{2\pi^2}}
\le c_4^2\|u\|_X^2.
\end{aligned}
\]
The estimates then follow from density of finite harmonic sums.
\end{proof}

The space $X$ consists of $H^2$-functions, so its norm controls the $L^2$-norm of first derivatives.
We record the following explicit bound for this:

\begin{proposition} \label{prop:sobolevderivatives}
For $u\in X$,
\[
\|d u\|_{L^2}\le c_5 \|u\|_X
\]
where
\[
c_5=\tfrac{2\sqrt{2}}{9}.
\]
\end{proposition}
\begin{proof} Let $u=\sum_k a_k\Psi_k$. Then

$$\|du\|_{L^2}^2=\sum_k\frac{l_k(l_k+2)}{(l_k+1)^4}\,(l_k+1)^4a_k^2\le\Big(\max_{l\ge2}\frac{l(l+2)}{(l+1)^4}\Big)\sum_k(l_k+1)^4a_k^2=\frac8{81}\|u\|_X^2 .$$
\end{proof}

\subsection{The harmonic map equation}
\label{subsection:harmonic-map-equation}

A smooth map $u:S^3\to S^2$ is harmonic if it is a critical point of the Dirichlet energy
\[
E(u)=\frac12\int_{S^3}|du|^2\,d\mathrm{vol}.
\]
This is the case if and only if it satisfies the elliptic equation
\begin{align}
\label{eqn:harmonic-map-eqn}
\tau(u)=\Delta u+|du|^2u=0,
\end{align}
see, with the opposite Laplacian convention, \cite[(1.12)--(1.13)]{FeehanMaridakis}.

Throughout the paper, we will solve the following equation, which is particularly amenable for numerically verified proofs.
For fixed $\bar W \in X$ with $\bar W$ nowhere zero define
\begin{align}
F:X&\longrightarrow Y,\nonumber\\
v&\longmapsto
\frac12\bar W\times\Big(W\times\big(3|W|^2\Delta W-\Delta(|W|^2W)\big)\Big)
+|\bar W|^2\big(\bar W\cdot(1-\Delta)v\big)\bar W,\label{eqn:F-equation}
\end{align}
where $W=\bar W+v$.
One checks that $F$ indeed maps equivariant functions to equivariant functions and, using the Sobolev embedding theorems, that it maps $H^2$-fields to $L^2$-fields.

\begin{proposition}\label{prop:zeros}
	If $v \in X$ satisfies $\|v\|_{C^0} < \min_{S^3} |\bar W|$ and $F(v)=0$, then 
	\[
	u=\frac{\bar W+v}{|\bar W+v|}
	\]
	is a smooth $G$-equivariant harmonic map.
\end{proposition}

\begin{proof}
Our assumption $\|v\|_{C^0} < \min_{S^3} |\bar W|$ ensures $W= \bar W + v$ is non-vanishing, so $u$ is well-defined. We have
\begin{align}
\label{eqn:normalization-identities}
u\times du=\frac{W\times dW}{|W|^2},\qquad
\operatorname{div}(u\times du)=u\times\Delta u,\qquad
\operatorname{div}(W\times dW)=W\times\Delta W.
\end{align}
This gives
\begin{align*}
|W|^4u\times\Delta u
&=|W|^4\operatorname{div}
  \left(\frac{W\times dW}{|W|^2}\right)\\
&=W\times\left(
  |W|^2\Delta W-\langle d(|W|^2),dW\rangle
  \right)\\
&=\frac12W\times
  \left(3|W|^2\Delta W-\Delta(|W|^2W)\right),
\end{align*}
where in the first step we used \eqref{eqn:normalization-identities};
in the second step we used the product rule,
$d(|W|^{-2})=-|W|^{-4}d(|W|^2)$,
and \eqref{eqn:normalization-identities};
and in the third step we used the product rule for
$\Delta(|W|^2W)$ and $W\times W=0$.

Plugging this into the definition of $F$ gives
\[
0=F(v)=
|W|^4\bar W\times(u\times\Delta u)
+|\bar W|^2\bigl(\bar W\cdot(1-\Delta)v\bigr)\bar W.
\]
The two summands are orthogonal, so $F(v)=0$ implies
$u\times\Delta u=0$.
Here we used that $\bar W\cdot u>0$ and
$u\times\Delta u$ is perpendicular to $u$.
Differentiating $|u|^2=1$ gives
\[
0=\frac12\Delta(|u|^2)u
=u\times(u\times\Delta u)+\Delta u+|du|^2u
=\Delta u+|du|^2u,
\]
which is \eqref{eqn:harmonic-map-eqn}.
Since $du\in L^6$, standard elliptic theory gives that $u$ is smooth.
\end{proof}

\begin{remark} \label{rem:gaugefix} As the proof of Proposition \ref{prop:zeros} shows, the two terms in \eqref{eqn:F-equation} vanish independently. The first term comes from the harmonic map equation \eqref{eqn:harmonic-map-eqn}, re-written for $u = W|W|^{-1}$, while the second term is a gauge-fixing, which deals with the freedom to rescale $W$ by a function to give the same harmonic map $u$. In this case, the second term in $\eqref{eqn:F-equation}$ fixes the component of $W = \bar{W} +v$ along $\bar{W}$ while the first term fixes the component of $W$ orthogonal to $\bar{W}$.  
\end{remark}
\subsection{The Hopf invariant}
\label{subsection:hopf-invariant}
\begin{definition}
Let $\omega$ be an area form on $S^2$ normalized by $\int_{S^2}\omega=1$.
Since $H^2_{dR}(S^3) =0$, for any smooth map $f:S^3\to S^2$, one can choose $\alpha\in\Omega^1(S^3)$ with $d\alpha=f^*\omega$. We define the (absolute) \textit{Hopf invariant} by
\[
H(f)=\left| \int_{S^3}\alpha\wedge d\alpha \right| \in\Z_{\geq 0}.
\]
This is independent of $\alpha$, invariant under homotopy, and extends to continuous maps by smooth approximation.
\end{definition}

One can compute the Hopf invariant in many cases of interest using the following proposition, proved in \cite[Theorem 2]{Hopf}: 
\begin{proposition} \label{prop:hopfinvariant} Suppose  $F:S^3\to S^3$, $f:S^3\to S^2$ and $G:S^2\to S^2$ are smooth, then
$$H(G\circ f\circ F)=(\deg G)^2\,|\deg F|\,H(f).$$
\end{proposition}
This is particularly useful  combined with the following map: 
\begin{definition} \label{def:hopffibration} Identify $\R^4=\C^2$ and $\R^3=\C\times\R$. For $z,w \in \C$, define the (extended) \textit{Hopf fibration} $\pi:\R^4 \rightarrow \R^3$ as the map
  \begin{align}
\pi(z,w)=(w^2-\bar z^2,2\operatorname{Re}(zw))   
\end{align}  
This restricts to a map  $\pi:S^3 \rightarrow S^2$ on the sphere $|z|^2+|w|^2=1$. 
\end{definition}
As a map $\pi: S^3 \rightarrow S^2$, \cite[Theorem 2]{Hopf} shows that the Hopf fibration satisfies $H(\pi) =1$. This has the following consequence 
\begin{lemma} \label{lem:hopfdegree} Suppose $F: S^3 \rightarrow \R^4$ with $|F|$ non-zero. Let $V = \pi(F)$, where $\pi$ is the extended Hopf fibration. Then $\tfrac{V}{|V|}$ is a map $S^3 \rightarrow S^2$, with 
\begin{align*}
H\left(\tfrac{V}{|V|}\right) = \left|\deg\left(\tfrac{F}{|F|}\right)\right|.
\end{align*}
\end{lemma}
\begin{proof}
This uses two key properties of the extended Hopf map $\pi$, which are easy to verify. These are that $|\pi(z,w)| = |z|^2+ |w|^2$, and that $\pi(\lambda z,\lambda w) = \lambda^2\pi(z,w)$ for any $\lambda>0$. In particular 
\begin{align*}
\frac{V}{|V|} = \frac{V}{|F|^2} = \pi \left(\frac{F}{|F|} \right).
\end{align*}
We now apply \cref{prop:hopfinvariant} and conclude.
\end{proof}

Later on, we will show that our approximate harmonic map is homotopic to a reference map with a known Hopf invariant, which will allow us to compute the Hopf invariant of the approximate solution.

The reference map we will use is the following $G_q$-equivariant map: 
\begin{align} \label{eq:referencemap}
    V_q(z,w) := \big(w^2-\tfrac{1}{4}\bar z^{2q},\,\Re(z^qw)\big)
\end{align}
Some properties of this reference map follow by basic calculation:
 \begin{lemma}  \label{lem:homotopysup} The map $V_q$ in \eqref{eq:referencemap} satisfies $|V_q| =\tfrac{1}{4}|z|^{2q}+|w|^2$. In particular, restricted to sphere $|z|^2+|w|^2=1$, we have the inequality $0<|V_q|\leq1$. 
 \end{lemma}
\begin{proof} The lower bound is trivial. On the sphere $|z|^2 \leq |z|^2 + |w|^2 =1$ so $\tfrac{1}{4}|z|^{2q}+|w|^2 \leq |z|^2+|w|^2 =1$, which proves the upper bound.   \end{proof}
Since $V_q$ is non-vanishing, we may consider $\tfrac{V_q}{|V_q|}$ as a map $S^3 \rightarrow S^2$, and use \cref{lem:hopfdegree} to compute its Hopf invariant. 
\begin{proposition} 
The Hopf invariant of $\tfrac{V_q}{|V_q|}$ satisfies $H\left(\tfrac{V_q}{|V_q|}\right) =q$.
\end{proposition} 
\begin{proof} Up to scaling $z$ by a constant (which does not change the homotopy class), $V_q$ is precisely $\pi (z^q,w)$. The map 
\begin{align*}
(z,w) \mapsto (z^q,w)    
\end{align*}
    suitably normalised, is a $q$-fold (branched) covering map on $S^3$, since on the regular locus $z \neq 0$, the pre-image of a point $(z,w)$ is generated by the $q^{th}$ root of unity. Hence, the degree of this map is precisely $q$. The proposition now follows by \cref{lem:hopfdegree}.
\end{proof}
\begin{remark}
The $z$-scaling factor in $V_q$ was chosen in order to ensure the numerical homotopy lower bound in \cref{prop:homotopy} below was sufficiently large, see also \cref{prop:approximate-bounds}, but its exact value is not important.  
\end{remark}
The following proposition reduces the construction of a homotopy to exhibiting a suitable vector field:

\begin{proposition}\label{prop:homotopy}Suppose that $W_0,W_1 :S^3\to\R^3$ are continuous and
\[
W_0\cdot W_1>0
\qquad\text{on }S^3.
\]
Then $W_0$, $W_1$ are nowhere zero, and $\frac{W_0}{|W_0|}$ is homotopic to $\tfrac{W_1}{|W_1|}$. In particular, $H(\frac{W_0}{|W_0|}) = H(\frac{W_1}{|W_1|})$. 
\end{proposition}
\begin{proof}
Firstly, $|W_0|, |W_1|$ must be  non-vanishing by Cauchy-Schwarz, as is $|(1-t)W_0 + t W_1|$ for $t\in [0,1]$. Thus, the following defines a homotopy between $W_0/|W_0|$ and $W_1/|W_1|$ as maps $S^3 \rightarrow S^2$:
\begin{align*}
 W_t =\tfrac{(1-t)W_0 + tW_1}{|(1-t)W_0 + tW_1|}.
\end{align*}
\end{proof}  

\subsection{Energy of a harmonic map \texorpdfstring{$R\circ\pi\circ A$}{R o π o A}} \label{subsection:energy}
In order to compare the energies of our harmonic maps to the energies of the known harmonic maps, we compute the energy of a harmonic map $u:S^3 \rightarrow S^2$ of the form $u = R \circ \pi \circ A$, where $A$ is an isometry of $S^3$, $\pi$ is the Hopf fibration, and $R:S^2 \rightarrow S^2$ is a weakly conformal map. Unfortunately, the authors could not find an explicit reference for this in the literature. To carry out this computation, we prove a short lemma:

\begin{lemma} \label{lem:hopfintegration} For any  $f \in L^1(S^2, \R)$, $\int_{S^3} \pi^*f = \tfrac{\pi}{2} \int_{S^2} f$.
\end{lemma}
\begin{proof} Let $s: S^2 \rightarrow S^2\left(\tfrac{1}{2}\right)$ be the dilation map $s(\zeta,t)=\left( \tfrac{1}{2}\zeta, \tfrac{1}{2}t\right)$ that sends the two sphere $S^2$ of radius $1$ to the two-sphere $S^2\left(\tfrac{1}{2}\right)$ of radius $\tfrac{1}{2}$. 

We use a standard fact about the Hopf fibration: the map $s\circ\pi: S^3 \rightarrow S^2\left(\tfrac{1}{2}\right)$ is a Riemannian submersion with fibres of constant length $2\pi$. In this case, the co-area formula for Riemannian submersions gives the required result:
\begin{align*}
\int_{S^3} \pi^*f = \int_{S^3} (s \circ\pi)^*(\left(s^{-1}\right)^*f )  = 2\pi \int_{S^2\left(\tfrac{1}{2}\right)} \left(s^{-1}\right)^*f = \tfrac{\pi}{2} \int_{S^2} f,
\end{align*}
where we have used the co-area formula in the second step, and in the third step, we have used that the volume form on $S^2\left(\tfrac{1}{2}\right)$ is $\tfrac{1}{4}\omega$, where $\omega$ is the volume form on $S^2$.
\end{proof}
\begin{proposition} \label{prop:hopfenergy} Let $u = R \circ \pi \circ A$ be a map $S^3 \rightarrow S^2$, where $A$ is an isometry of $S^3$, $\pi$ is the Hopf fibration, and $R:S^2 \rightarrow S^2$ is a weakly conformal map. Then 
$$E (u) = \tfrac{1}{2} \int_{S^3} |du|^2 = 8\pi^2 |\deg (R)| $$
\end{proposition}
\begin{proof} First, since $A$ is an isometry $E(R \circ \pi \circ A) = E(R \circ \pi)$, so it suffices to compute the energy of a map of the form $u = R \circ \pi $. Now, since $s \circ \pi$ is a Riemannian submersion $S^3 \rightarrow S^2\left(\tfrac{1}{2}\right)$, it follows that $d\pi$ acts on horizontal vector-fields of the Hopf fibration as multiplication by $2$. Therefore, for all $x \in S^3$, $|du_x|^2 = |d R_{\pi(x)}\circ d\pi_x|^2 = 4|dR_{\pi(x)}|^2$. Applying \cref{lem:hopfintegration}, gives 
\begin{align*}
E (u) = 2\int_{S^3} \pi^* |dR|^2 =  \pi \int_{S^2}|dR|^2 = 2 \pi E(R)  
\end{align*}
Now, a weakly conformal map $R : S^2 \rightarrow S^2$ is automatically harmonic, we have that $E(R) = 4 \pi | \deg(R)|$ by \cite[p.130]{EellsSampson}, and we are done. 
\end{proof}
\subsection{Index of a harmonic map}\label{subsection:index}

The second variation of the Dirichlet energy functional, at a critical point $u:S^3 \rightarrow S^2$ with $\tau(u) = 0$ is given by the quadratic form: 

\begin{align} \label{eq:quadraticform0}
Q_u (\xi) = \int_{S^3} |d\xi|^2 - |du|^2|\xi|^2  
\end{align}
where $\xi \cdot u = 0$ and $d \xi$ is the component-wise derivative of $\xi$ in  $\mathbb{R}^3$. This quadratic form is defined on sections $\Gamma_{H^1}(u^*TS^2) = \{\xi\in H^1(S^3,\mathbb R^3) : \xi\cdot u = 0\}$.

In this section, we will work with the following quadratic form, which is better adapted to our set up. For fixed $W \in X$ with $ W$ nowhere zero, and $\xi\in H^1(S^3,\mathbb R^3)$, define $Q_W: H^1(S^3,\mathbb R^3) \rightarrow \mathbb R$ as
\begin{align} \label{eq:quadraticform} 
Q_W(\xi):=\int_{S^3}\Big|d\xi-\frac{(W\times dW)\times\xi}{|W|^2}\Big|^2-\frac{|(W\times dW)\cdot\xi|^2}{|W|^4}+\frac{(W\cdot\xi)^2}{|W|^2}.\end{align}
\begin{proposition}\label{prop:index-equivalence} Let $Q_u$, $Q_W$, be the quadratic forms in \eqref{eq:quadraticform0}, \eqref{eq:quadraticform} acting on  $\Gamma_{H^1}(u^*TS^2) = \{\xi\in H^1(S^3,\mathbb R^3) : \xi\cdot u = 0\}$ and $H^1(S^3,\mathbb R^3)$ respectively, and suppose $u = \tfrac{W}{|W|}$. Then: 
\begin{align*}
    \mathrm{Ind} Q_u = \mathrm{Ind} Q_W & &\ker Q_u = \ker Q_W
\end{align*}
\end{proposition}
\begin{proof} Let $P_u = I - u u^T$ be the projection to the orthogonal complement of $u$. We defined $Q_W$ in \eqref{eq:quadraticform} such that
\begin{align}
Q_W (\xi) = Q_u(P_u\xi) + \int_{S^3} |d (u\cdot\xi)|^2 + \left(u\cdot\xi\right)^2 .
\end{align}
  From this, one can see that $Q_u$ and $Q_W$ agree on $P_u(\xi)$, the component of $\xi$ orthogonal to $u$. They differ only on the component parallel to $u$, but by a positive-definite form. The parallel summand is positive definite, so it contributes no index or kernel.
\end{proof}
The point of doing this is so that we can pass to a quadratic form over arbitrary maps $\xi$ into $\mathbb R^3$, just like we did when we gauge-fixed the components in the $\bar{W}$-directions in equation \eqref{eqn:F-equation}, see \cref{rem:gaugefix}. 

\subsubsection{Estimating the index}

Later on we will use an approximate harmonic map $\bar W$ and perturb it to a genuine harmonic map $W$.
To ensure that negative directions for $Q_{\bar W}$ remain negative for $Q_W$, we estimate the difference between these quadratic forms.

\begin{proposition} \label{prop:quadraticestimate} Let $W,\bar W\in X$, $r,\Cr{const:energyparameter}>0$ with $\|\bar W-W\|_X\le r$ and $\min_{S^3}|\bar W|-c_\infty r
\ge \Cl{const:energyparameter}>0$. Then, for $\xi\in H^1(S^3,\mathbb R^3)$, $|Q_W(\xi)-Q_{\bar W}(\xi)|\le\delta\|\xi\|_{H^1}^2$ with
$$
\delta = r\Big[2\beta\Big(\frac{3c_\infty c_4(\|\bar W\|_X+r)}{\Cr{const:energyparameter}^2}+\frac{c_4}{\Cr{const:energyparameter}}\Big)+\frac{12\beta^2c_\infty c_4^2(\|\bar W\|_X+r)^2}{\Cr{const:energyparameter}^3}+\frac{4\beta^2c_4^2(\|\bar W\|_X+r)}{\Cr{const:energyparameter}^2}+\frac{4 c_\infty}{\Cr{const:energyparameter}}\Big].
$$
where $\beta$, $c_4$, $c_\infty$ are the constants in Proposition \ref{prop:sobolev}. 
\end{proposition}
\begin{proof}
 Write $W=\bar W+v$ with $\|v\|_X\le r$, and define 
$$A:=\frac{W\times dW}{|W|^2},\qquad u:=\frac{W}{|W|},$$
with $\bar{A}$, $\bar{u}$ defined similarly for $\bar W$. We need to bound
\begin{equation}
\begin{aligned} \label{eq:quadraticbound}
Q_W(\xi)-Q_{\bar W}(\xi)&=\int-2\langle d\xi,(A-\bar A)\times\xi\rangle+\big(|A\times\xi|^2-|\bar A\times\xi|^2\big) \\&- \big(|A\cdot\xi|^2-|\bar A\cdot\xi|^2\big)+\big((u\cdot\xi)^2-(\bar u\cdot\xi)^2\big).    
\end{aligned}
\end{equation}
in terms of $\|\xi\|_{H^1}$. 

We can factorise the last three terms in \eqref{eq:quadraticbound} as the \textit{difference of squares}: e.g. for the third term in the integrand of  \eqref{eq:quadraticbound},
\begin{align*}
 |A \cdot \xi|^2-|\bar A\cdot\xi|^2= \big((A-\bar A)\cdot\xi\big)\cdot\big((A+\bar A)\cdot\xi\big).  
\end{align*}
Using this factorisation, we see that the integrand in \eqref{eq:quadraticbound} is bounded in absolute value by 
\begin{align} \label{eqn:quadraticpointwise} 
2|A-\bar A|\,|\xi||d\xi|+2|A-\bar A|\big(|A|+|\bar A|\big)|\xi|^2+|u+\bar{u}||u-\bar u|\,|\xi|^2.    
\end{align}
In the following,  we define the following functions on $S^3$:
\begin{align} \label{eqn:quadraticterms}
T_1 :=  2|A-\bar A|,& &T_2:=2|A-\bar A|\big(|A|+|\bar A|\big)+|u+\bar{u}||u-\bar u|.
\end{align}
Our proof strategy is as follows:  we will apply H\"{o}lder's inequality with three exponents $\tfrac{1}{4} + \tfrac{1}{4}+ \tfrac{1}{2} = 1$ to the integral of \eqref{eqn:quadraticpointwise}:
\begin{align*}
\int_{S^3}T_1|\xi||d\xi|+T_2|\xi|^2\le\|T_1\|_{L^4}\|\xi\|_{L^4}\|d\xi\|_{L^2}+\|\xi\|_{L^4}^2\|T_2\|_{L^2}\le\big(\beta\|T_1\|_{L^4}+\beta^2\|T_2\|_{L^2}\big)\|\xi\|_{H^1}^2
\end{align*}
where we have used Proposition \ref{prop:sobolev} in the second inequality. Since the integral of \eqref{eqn:quadraticpointwise} bounds \eqref{eq:quadraticbound}, this gives us the bound $|Q_W(\xi)-Q_{\bar W}(\xi)|\le\delta\|\xi\|_{H^1}^2$ in the statement of the proposition if we can bound $\|T_1\|_{L^4}$, and $\|T_2\|_{L^2}$. 

To bound $\|T_1\|_{L^4}$, and $\|T_2\|_{L^2}$, we first give a pointwise bound on $T_1$, $T_2$, which we will combine with the global bound $|v| \le c_\infty r$ (coming from \cref{prop:sobolev}) where appropriate.
\begin{enumerate}
    \item Bounds on $u$. Since $|u|=|\bar u | = 1$ by definition, $|u+\bar u| \leq 2$. Using $|a/|a|-b/|b||\le2|a-b|/|a|$ for nonzero $a,b$, we have that  $$|u-\bar u| \leq \frac{2|v|}{\Cr{const:energyparameter}} \le \frac{2c_\infty r}{\Cr{const:energyparameter}}. $$ 
    \item Bounds on $A$. Since $|W \times d W| \leq |W||dW|$, we have that 
    $$|A|\leq \frac{|d W|}{\Cr{const:energyparameter}}, \qquad |A|+|\bar A| \leq 2\frac{|d\bar W|+|dv|}{\Cr{const:energyparameter}}.$$
    For $A-\bar A$, let $W_t := \bar{W} + t v$ and differentiate $A_t=W_t\times dW_t/|W_t|^2$ in $t$:

$$\dot A_t=\frac{v\times dW_t+W_t\times dv}{|W_t|^2}-\frac{2(W_t\cdot v)\,W_t\times dW_t}{|W_t|^4}.$$

Since $|d W_t| \leq |d\bar W|+|dv|$, and $|W_t| \geq \Cr{const:energyparameter}$, so $|\dot A_t|\le\dfrac{3 |v| \left( |d\bar W|+|dv|\right)}{\Cr{const:energyparameter}^2}+\dfrac {|dv|}{\Cr{const:energyparameter}}$. Integrating over $t\in[0,1]$ gives

$$|A-\bar A|\le\dfrac{3 |v| \left( |d\bar W|+|dv|\right)}{\Cr{const:energyparameter}^2}+\dfrac {|dv|}{\Cr{const:energyparameter}}\leq \dfrac{3 r c_\infty \left( |d\bar W|+|dv|\right)}{\Cr{const:energyparameter}^2}+\dfrac {|dv|}{\Cr{const:energyparameter}}$$
\end{enumerate}

Equipped with these pointwise estimates, we can finally bound $\|T_1\|_{L^4}$, since, by Proposition \ref{prop:sobolev}: 
\begin{align}
\|dv\|_{L^4}\le c_4r,& &\||d \bar W| + |dv|\|_{L^4}\le c_4(\|\bar W\|_X+r).
\end{align}
Similarly for $\|T_2\|_{L^2}$, we use \cref{prop:sobolev} to conclude: 
\begin{align}
\||d v|(|d \bar W| + |dv|)\|_{L^2}\le\|dv\|_{L^4}\||d \bar W| + |dv|\|_{L^4}  \leq  r c_4^2  (\|\bar W\|_X+r)
\end{align}  
where in the first step, we use H\"{o}lder's inequality with exponents $\tfrac{1}{2} = \tfrac{1}{4}+ \tfrac{1}{4}$. Combining these estimates shows that we can indeed take $\delta$ as stated. 
\end{proof}
We leverage \cref{prop:quadraticestimate} in the following way. Let $Q(\xi,\eta)$ be the symmetric bilinear form built out of the quadratic form $Q$, i.e. 
$$Q(\xi,\eta):=\tfrac12\left[Q(\xi+\eta)-Q(\xi)-Q(\eta)\right].$$
Let $\xi_1,\ldots,\xi_k\in H^1(S^3,\mathbb R^3)$, and let $G,H$ be the $k\times k$ symmetric matrices with entries: 

$$G_{ab}=Q_{\bar W}(\xi_a,\xi_b),\qquad H_{ab}=\int_{S^3}\xi_a\cdot\xi_b+\langle d\xi_a,d\xi_b\rangle .$$
\begin{lemma}\label{lem:negdef}
Under the assumptions of \cref{prop:quadraticestimate}, let $\xi_1,\dots,\xi_k\in H^1(S^3,\mathbb R^3)$ and let $G,H$ be as above. If $G+\delta H$ is negative definite, then $Q_W$ is negative definite on their span. In particular $\operatorname{Ind}Q_W\ge k$.
\end{lemma}

\begin{proof}
Let $c\in\mathbb R^k\setminus\{0\}$ and $\xi=\sum_a c_a\xi_a$. By bilinearity,
$Q_{\bar W}(\xi)=c^TGc$ and $\|\xi\|_{H^1}^2=c^THc$, so \cref{prop:quadraticestimate} gives
\[
  Q_W(\xi)\le Q_{\bar W}(\xi)+\delta\|\xi\|_{H^1}^2=c^T(G+\delta H)c<0 .
\]
In particular, the $\xi_a$ are linearly independent.
\end{proof}

\subsection{The fixed point theorem}
\label{subsection:fixed-point}

Our existence proof proceeds by constructing an explicit approximate solution, and then proving by means of a fixed-point argument that a nearby exact solution exists.
This standard approach commonly uses the Newton-Kantorovich theorem as a sufficient criterion for the existence of the true solution, see e.g. \cite[Theorem~1]{SekineTakayasuOishi}.
For simpler statements, we use the following straightforward consequence of the Banach fixed point theorem:

\begin{theorem}\label{thm:NK}
Let $F:X\to Y$ be a continuously differentiable map between Banach spaces, and suppose
$\mathscr L=DF(0)$ is invertible. Set $\mathscr N(v)=F(v)-F(0)-\mathscr L v$.
Suppose that, for some $r>0$,
\[
\|\mathscr L^{-1}\|_{Y\to X}\le M,\qquad \|F(0)\|_Y\le\varepsilon,\qquad
\|\mathscr N(v)-\mathscr N(w)\|_Y\le k\|v-w\|_X
\]
whenever $\|v\|_X,\|w\|_X\le r$.
If $Mk<1$ and $M(\varepsilon+kr)\le r$, then $F$ has a unique zero in the closed ball of radius $r$ about $0$.
\end{theorem}

\section{Estimates}
\label{section:estimates}

In \cref{thm:main} we claim existence of five different harmonic maps.
Accordingly, the proof uses five different candidate maps and we compute estimates for each of the candidates.
Throughout, we abide by the convention of using the same symbols as placeholders that can be filled with the objects pertaining to each of the candidates.
For example, we will write $\bar W$ for an approximate harmonic map, and it is understood that this is a placeholder for five approximate solutions.

\subsection{Computing the approximate solution}
\label{subsection:computing-the-approx-sol}

Using the basis from \cref{prop:harmonic-bases}, we aim to find $\bar W=\sum_{l_k\le D}a_k\Psi_k$ solving
\[
  \big\langle\Psi_k,F(0)+10(|\bar W|^2-1)\bar W\big\rangle_{L^2}=0,
  \qquad l_k\le D.
\]
Here, the two terms are pointwise orthogonal,
since $\bar W\cdot F(0)=0$, and the second penalises deviations
from unit length.

The integrands in this system and its Jacobian are $G_q$-invariant polynomials of degree at most $6D+2$.
We evaluate these integrals using the exact quadrature from \cref{prop:orbit-quadrature}.
We then approximately solve the system using damped Newton iteration.
Quadrature evaluation and Newton iteration use non-rigorous double precision arithmetic.
We point out that there is no theoretical guarantee that the iteration converges to a (truncated) true solution, or even that the non-rigorous arithmetic gives results anywhere near the rigorous values.
All rigour comes from the following a posteriori analysis.

The degrees and leading coefficients are given in
\cref{master-table,tab:coefficients}.
The stored floating-point coefficients are treated as exact
rational numbers in the subsequent verification.

\subsection{Bounds for the approximate solution}

As our approximate solution, take the polynomial field $\bar W:S^3\to\R^3$ as given in the appendix, for a given $q$. Let $V = V_q$ be the homotopy reference map \eqref{eq:referencemap} in \S \ref{subsection:hopf-invariant}.

\begin{proposition}\label{prop:approximate-bounds}
For $r$ as in \cref{master-table}, the approximate solution $\bar W$
satisfies
\begin{align*}
\|F(0)\|_Y
&\le\Cl{const:residual},\\
\Cl{const:energy-lower}
&\le E\left(\frac{\bar W}{|\bar W|}\right)
\le\Cl{const:energy-upper},\\
\min_{S^3}\bar W\cdot V_q
&\ge\Cl{const:homotopy}>0,\\
\|\bar W\|_X
&\le\Cl{const:Wbar-X},\\
\|\bar W\|_{C^0}
&\le\Cl{const:Wbar-C0},\\
\big\||\bar W|^4-1\big\|_{C^0}
&\le\eta<1,\\
\min_{S^3}|\bar W|-c_\infty r
&\ge\Cr{const:energyparameter}>0,
\end{align*}
where the constants are given in \cref{master-table}.
\end{proposition}

\begin{proof}
    Integrals of $G_q$-invariant polynomials are computed using \cref{prop:orbit-quadrature}.
    Suprema and infima are bounded using the Lipschitz method (see e.g. \cite[p.~312, equation~(4)]{PaulaviciusZilinskas2009} for a description of this standard approach):
    evaluate on a finite grid and control values between grid points using derivative bounds computed using a Chebyshev--Fourier expansion.
    All computations are carried out using interval arithmetic.
    This suffices to compute $\Cr{const:residual},\Cr{const:homotopy},\Cr{const:Wbar-C0},\eta,\Cr{const:energyparameter}$.

    The constant $\Cr{const:Wbar-X}$ is computed using \cref{eqn:X-norm-calculation}.

    For $\Cr{const:energy-lower},\Cr{const:energy-upper}$, use $N=|\bar W\times d\bar W|^2$ and 
    \[
    \begin{aligned}
    \left|
    E\!\left(\frac{\bar W}{|\bar W|}\right)
    -\frac12\int_{S^3}N(2-|\bar W|^4)
    \right|
    &=\frac12\int_{S^3}
    \frac{N(|\bar W|^4-1)^2}{|\bar W|^4}
    \le
    \frac{\eta^2}{2(1-\eta)}\int_{S^3}N,
    \end{aligned}
    \]
    where the final integral is again computed using \cref{prop:orbit-quadrature}.
\end{proof}

\begin{proposition}\label{prop:indexbound}
For $q=2$, there exist nine linearly independent fields $\xi_1,\ldots,\xi_9\in H^1(S^3,\mathbb R^3)$, orthogonal for both $Q_{\bar W}$ and $\langle\cdot,\cdot\rangle_{H^1}$, such that
\[
\max_{1\le a\le9}\frac{Q_{\bar W}(\xi_a)}{\|\xi_a\|_{H^1}^2}\le\lambda^*<0,
\]
where $\lambda^*$ is given in \cref{master-table}.
\end{proposition}

\begin{proof}
By abuse of notation, we denote the Hermitian extension of both quadratic forms $Q_{\bar W}$ and $\langle\cdot,\cdot\rangle_{H^1}$ to complex fields by the same symbols. Let $f_1$, $f_2$, $f_3$, $f_4$ and $f_5$ be the nonzero elements of $H^1(S^3,\mathbb C^3)$ given in {\color{blue}\cref{tab:index-fields}}, with $f_1$ real. Take the nine real maps
\[
\xi_1=f_1,\qquad \xi_{2a-2}=\operatorname{Re}f_a,\quad \xi_{2a-1}=\operatorname{Im}f_a\qquad(a=2,\dots,5).
\]
By appendix \ref{appendix:index-fields}, these are non-zero and pairwise orthogonal for both $Q_{\bar{W}}$ and $\langle\cdot,\cdot\rangle_{H^1}$, hence linearly independent. For $a=2,\dots,5$, appendix \ref{appendix:index-fields} gives
\[
\frac{Q_{\bar W}(f_a)}{\|f_a\|_{H^1}^2}
=
\frac{Q_{\bar W}(\Re f_a)}{\|\Re f_a\|_{H^1}^2}
=
\frac{Q_{\bar W}(\Im f_a)}{\|\Im f_a\|_{H^1}^2}.
\]
Since $f_1$ is real, it therefore suffices to enclose the five
Rayleigh quotients of $f_1,\dots,f_5$.

We enclose the results of these computations in \cref{master-table}, which are carried out as follows. 
All integrands are sums of circle-invariant polynomials times $1$, $|\bar W|^{-2}$ or $|\bar W|^{-4}$.
We replace these factors by $1$, which gives an error controlled by $\sup_{S^3}\big||\bar W|^2-1\big|$, and integrate the polynomials exactly using \cref{prop:orbit-quadrature}. 
The resulting enclosure gives the stated bound.
\end{proof}

\subsection{The linear estimate}\label{subsection:linear-estimate}

The linearisation of $F$ at $v=0$ is given by
\begin{align}
\mathscr L h=DF(0)h
={}&\frac12\bar W\times\Big[
h\times\big(3|\bar W|^2\Delta\bar W-\Delta(|\bar W|^2\bar W)\big)\nonumber\\
&\qquad+\bar W\times\big(6(\bar W\cdot h)\Delta\bar W+3|\bar W|^2\Delta h
-\Delta(2(\bar W\cdot h)\bar W+|\bar W|^2h)\big)\Big]\nonumber\\
&+|\bar W|^2\big(\bar W\cdot(1-\Delta)h\big)\bar W.
\label{eq:linearisation}
\end{align}

\subsubsection{An \texorpdfstring{$L^2$}{L²}-inverse bound}

Fix $l\in\Z_{\ge0}$ and define the spaces of \emph{low modes} and \emph{high modes} as
\[
\mathcal H_{\le l}:=\operatorname{span}\{\Psi_k:l_k\le l\},\qquad
\mathcal H_{>l}:=\operatorname{span}\{\Psi_k:l_k>l\}.
\]
The $H^2$-completion of $\mathcal{H}_{\leq l} \oplus \mathcal{H}_{>l}$ is $X$, and according to this splitting we write
\begin{align*}
\mathscr L =
\begin{pmatrix}
A&B\\
C&D
\end{pmatrix}:\mathcal{H}_{\leq l} \oplus \mathcal{H}_{>l} & \rightarrow \mathcal{H}_{\leq l} \oplus \mathcal{H}_{>l}.
\end{align*}

\textbf{Lower bound for $D$.} Write $\mathscr R=\mathscr L+\Delta$.
The second order term of $\mathscr R$ is $(|\bar W|^4-1)(-\Delta)$.
In a local orthonormal frame, write its lower order part as $\sum_i M_i\nabla_i+M_0$,
with coefficients obtained by expanding \eqref{eq:linearisation}.
Choose $b_1>0$ and $b_0>0$ such that
\begin{align}
\label{eqn:remainder-estimate}
\sup_{S^3}\sum_i\|M_i\|_{\mathrm F}^2\le b_1^2,\qquad
\sup_{S^3}\|M_0\|_{\mathrm F}^2\le b_0^2,
\end{align}
where $\|\cdot\|_{\mathrm F}$ is the Frobenius norm.
These expressions are independent of the frame. Hence
\[
\|\mathscr Rh\|_{L^2}
\le\eta\|\Delta h\|_{L^2}
+b_1\|\nabla h\|_{L^2}+b_0\|h\|_{L^2},
\]
where $\eta$ is as in \cref{prop:approximate-bounds}.
The bounds for $b_1,b_0$ are recorded in \cref{master-table}.
On $\mathcal H_{>l}$, $-\Delta$ has eigenvalues at least $(l+1)(l+3)$.
Hence, for $h \in \mathcal{H}_{>l}$:
\[
\|h\|_{L^2}\le\frac{\|\Delta h\|_{L^2}}{(l+1)(l+3)},\qquad
\|\nabla h\|_{L^2}^2=-\langle h,\Delta h\rangle_{L^2}
\le\frac{\|\Delta h\|_{L^2}^2}{(l+1)(l+3)}.
\]
Thus,
\[
\|\mathscr Rh\|_{L^2}\le
\left(\eta+\frac{b_1}{\sqrt{(l+1)(l+3)}}+\frac{b_0}{(l+1)(l+3)}\right)\|\Delta h\|_{L^2}.
\]
Put $\gamma:=(1-\eta)(l+1)(l+3)-b_1\sqrt{(l+1)(l+3)}-b_0$ and suppose $\gamma>0$.
Writing $\Pi$ for the $L^2$-orthogonal projection onto the high modes, we get
\begin{align*}
\|Dh\|_{L^2}
&=
\|\Pi\mathscr Lh\|_{L^2}
\\
&\geq
\|\Delta h\|_{L^2}-\|\Pi\mathscr Rh\|_{L^2}
\\
&\geq
\|\Delta h\|_{L^2}-\|\mathscr Rh\|_{L^2}
\\
&\geq
\frac{\gamma}{(l+1)(l+3)}\|\Delta h\|_{L^2}
\\
&\geq
\gamma\|h\|_{L^2},
\end{align*}
where in the third-to-last step we used that the orthogonal projection $\Pi$ has operator norm $\leq 1$.
On the high modes,
\begin{align}
\label{eqn:laplace-D-comparison}
D=\big(I+\Pi\mathscr R(-\Delta)^{-1}\big)(-\Delta),\qquad
\|\Pi\mathscr R(-\Delta)^{-1}\|_{L^2\to L^2}
\le1-\frac\gamma{(l+1)(l+3)}<1.
\end{align}
The Neumann series therefore proves that $D$ is invertible from high-mode $H^2$ to high-mode $L^2$, with
$\|D^{-1}\|_{L^2\to L^2}\le\gamma^{-1}$.

\textbf{Lower bound for $A$.}
We compute the entries
\[
A_{\alpha\beta}
=\int_{S^3}\Psi_\alpha\cdot\mathscr L\Psi_\beta\,d\mathrm{vol},
\qquad l_\alpha,l_\beta\le l,
\]
using rigorous quadrature and interval arithmetic and obtain an interval
matrix $\mathbf A\subset\mathbb R^{d\times d}$ containing $A$, where
$d=\dim\mathcal H_{\le l}$.

We compute an upper bound for $\|A^{-1}\|_2$ as follows:
write $m(\mathbf A)$ for its midpoint matrix and compute a non-rigorous
candidate for its inverse, $Z\approx m(\mathbf A)^{-1}$, using
floating-point arithmetic.
We define $\mathbf E:=I-Z\mathbf A$ and verify
\[
\|\mathbf E\|_2
\le\sqrt{\|\mathbf E\|_1\|\mathbf E\|_\infty}<1
\]
rigorously using interval arithmetic.
Here norms of interval matrices denote suprema over their elements.
By the Neumann series, every matrix $a\in\mathbf A$ is invertible with
\[
\|a^{-1}\|_2
\le\frac{\|Z\|_2}{1-\|\mathbf E\|_2}
\le
\frac{\sqrt{\|Z\|_1\|Z\|_\infty}}
{1-\sqrt{\|\mathbf E\|_1\|\mathbf E\|_\infty}},
\]
where the right-hand side is easily computed rigorously.
This is the same approach as \cite[Lemma~2.1, (2.7)]{RumpSingularValues2011}.

\textbf{Upper bounds for $B$ and $C$.}
Let $\mathscr{L}^*$ be the formal adjoint of $\mathscr{L}$.
Define the $d \times d$ matrices
\[
G_{\alpha \beta}
=
\langle \mathscr{L} \Psi_\alpha, \mathscr{L} \Psi_\beta \rangle_{L^2},
\quad
G_{\alpha \beta}^*
=
\langle \mathscr{L}^* \Psi_\alpha, \mathscr{L}^* \Psi_\beta \rangle_{L^2},
\quad
l_\alpha,l_\beta \leq l.
\]
Orthogonality of low mode and high mode spaces gives
\[
C^TC=G-A^TA,\qquad BB^T=G^*-AA^T.
\]
The matrix $A$ has been computed before, which gives $A^TA$ and $AA^T$.
The matrices $G$ and $G^*$ are readily computed using exact quadrature.
The squared norms \(\|C\|_2^2\) and \(\|B\|_2^2\) are the largest eigenvalues of \(G-A^TA\) and \(G^*-AA^T\), respectively.
We bound them using
\[
\lambda_{\max}(S)
\le
\min\left\{
\|S^{J}\|_{\mathrm F},
\sqrt{\|S^{J}\|_1\|S^{J}\|_\infty}
\right\}^{1/J},
\quad
J \in \mathbb{Z}_{\geq 1}
\]
for symmetric positive semidefinite $S$.
In the companion code we choose $J=32$.

\textbf{Combining the bounds using the Schur complement.}
All operator norms below are induced by $L^2$.
Using the same approximate inverse $Z$, we also bound
\[
\|A^{-1}B\|
\le \frac{\|ZB\|_2}{1-\|\mathbf E\|_2}
\quad
\text{ using }
\quad
\|ZB\|_2^2
=\|ZBB^TZ^T\|_2
=\|Z(G^*-AA^T)Z^T\|_2,
\]
where the finite-dimensional matrices on the right hand side have been computed before.
Let
$S:=A-BD^{-1}C:\mathcal H_{\le l}\to\mathcal H_{\le l}$
be the Schur complement of $D$.
For $u\in\mathcal H_{\le l}$,
\begin{align*}
\|Su\|_{L^2}
&\ge \frac{\|A^{-1}Su\|_{L^2}}{\|A^{-1}\|}\\
&=\frac{\|u-A^{-1}BD^{-1}Cu\|_{L^2}}{\|A^{-1}\|}\\
&\ge
\frac{1-\|A^{-1}B\|\|D^{-1}\|\|C\|}
{\|A^{-1}\|}\|u\|_{L^2}.
\end{align*}
Suppose the numerator is positive.
Since $\mathcal H_{\le l}$ is finite-dimensional, $S$ is invertible.
The block inverse formula gives
\begin{align}
\label{eqn:L-block-inverse}
\mathscr L^{-1}
&=
\begin{pmatrix}
S^{-1} & -S^{-1}BD^{-1}\\
-D^{-1}CS^{-1} & D^{-1}+D^{-1}CS^{-1}BD^{-1}
\end{pmatrix}
=
\begin{pmatrix}
0&0\\
0&D^{-1}
\end{pmatrix}
+
\begin{pmatrix}
I\\
-D^{-1}C
\end{pmatrix}
S^{-1}
\begin{pmatrix}
I&-BD^{-1}
\end{pmatrix}.
\end{align}
Taking norms, using the $L^2$-orthogonality of the splitting and
the preceding estimate, gives
\[
\|\mathscr L^{-1}\|_{L^2\to L^2}
\le
\|D^{-1}\|
+
\frac{
\|A^{-1}\|
\sqrt{
\left(1+\|B\|^2\|D^{-1}\|^2\right)
\left(1+\|C\|^2\|D^{-1}\|^2\right)
}
}{
1-\|A^{-1}B\|\|D^{-1}\|\|C\|
}.
\]

\subsubsection{A priori estimate}
For $h\in X$, recalling that $\eta<1$, we have
\begin{align*}
\|\Delta h\|_{L^2}
&\le \|\mathscr Lh\|_{L^2}+\|\mathscr Rh\|_{L^2}\\
&\le \frac{1}{1-\eta}
\left(\|\mathscr Lh\|_{L^2}
+b_1\|\nabla h\|_{L^2}+b_0\|h\|_{L^2}\right)\\
&\le \frac{1}{1-\eta}
\left(\|\mathscr Lh\|_{L^2}
+b_1\sqrt{\|h\|_{L^2}\|\Delta h\|_{L^2}}
+b_0\|h\|_{L^2}\right)\\
&\le \frac{\|\mathscr Lh\|_{L^2}}{1-\eta}
+\frac12\|\Delta h\|_{L^2}
+\left(\frac{b_0}{1-\eta}
+\frac{b_1^2}{2(1-\eta)^2}\right)\|h\|_{L^2}.
\end{align*}
Here the first step uses the definition of $\mathscr R$;
the second uses \eqref{eqn:remainder-estimate}, the preceding inequality,
and rearrangement; the third uses integration by parts and
Cauchy--Schwarz; and the fourth uses Young's inequality.

Thus, absorbing $\frac12\|\Delta h\|_{L^2}$, we obtain
\[
\|h\|_X
\le \|\Delta h\|_{L^2}+\|h\|_{L^2}
\le \frac{2}{1-\eta}\|\mathscr Lh\|_{L^2}
+\left(1+\frac{2b_0}{1-\eta}
+\frac{b_1^2}{(1-\eta)^2}\right)\|h\|_{L^2},
\]
where the first inequality follows from
$\|h\|_X=\|(1-\Delta)h\|_{L^2}$ and the triangle inequality.

\subsubsection{Injectivity estimate}
Combining the estimates from the preceding two subsubsections, we obtain:

\begin{proposition}
    \label{prop:injectivity-estimate}
    The operator $\mathscr{L}:X \rightarrow Y$ is invertible and satisfies $\|\mathscr{L}^{-1}\| \leq M$ for $M$ recorded in \cref{master-table}.
\end{proposition}

\begin{proof}
The operator $-\Delta: X \rightarrow Y$ is an isomorphism, also when restricted to high modes.
By \cref{eqn:laplace-D-comparison}, $D$ is an isomorphism.
Hence, using \cref{eqn:L-block-inverse}, one checks that $\mathscr{L}^{-1}$ is indeed a map $Y \rightarrow X$, proving that $\mathscr{L}$ is invertible.

The previous estimates give for $h\in X$,
\begin{align*}
\|h\|_X
&\le \frac{2}{1-\eta}\|\mathscr Lh\|_{L^2}
+\left(1+\frac{2b_0}{1-\eta}
+\frac{b_1^2}{(1-\eta)^2}\right)\|h\|_{L^2}\\
&\le
\left[
\frac{2}{1-\eta}
+\left(1+\frac{2b_0}{1-\eta}
+\frac{b_1^2}{(1-\eta)^2}\right)
\|\mathscr L^{-1}\|_{L^2\to L^2}
\right]
\|\mathscr Lh\|_{Y}
\\&\le
\left[
\frac{2}{1-\eta}
+\left(1+\frac{2b_0}{1-\eta}
+\frac{b_1^2}{(1-\eta)^2}\right)
\left(
\frac1\gamma+
\frac{
\|A^{-1}\|
\sqrt{\left(1+\frac{\|B\|^2}{\gamma^2}\right)
      \left(1+\frac{\|C\|^2}{\gamma^2}\right)}
}{
1-\frac{\|A^{-1}B\|\,\|C\|}{\gamma}
}
\right)
\right]
\|\mathscr Lh\|_Y.
\end{align*}
The constant $M$ is given by the expression in the square brackets above, in terms of our previous estimates.  
\end{proof}

\subsection{Higher order estimate}

It remains to prove an estimate for the higher order terms of \eqref{eqn:F-equation}.
In the following, write $\nabla_i$ for differentiation along a local orthonormal frame; all contractions below are frame-independent.
The higher order terms are then
\begin{align*}
\mathscr N:X&\longrightarrow Y,\\
v&\longmapsto F(v)-F(0)-\mathscr L v\\
&\phantom{\longmapsto{}}=6\mathcal T(\bar W,\bar W,v,v)+4\mathcal T(\bar W,v,v,v)+\mathcal T(v,v,v,v),\\
\mathcal T(a_1,a_2,a_3,a_4)
&=\frac1{24}\bar W\times\sum_{\sigma\in S_4}
\Big[(a_{\sigma(1)}\cdot a_{\sigma(2)})(a_{\sigma(3)}\times\Delta a_{\sigma(4)})
-2\sum_i(a_{\sigma(1)}\cdot\nabla_i a_{\sigma(2)})
(a_{\sigma(3)}\times\nabla_i a_{\sigma(4)})\Big],
\end{align*}
for $a=(a_1,a_2,a_3,a_4)\in X^{\oplus 4}$.

\begin{proposition}\label{prop:nonlinear}
For $v,w\in X$ with $\|v\|_X,\|w\|_X\le r$, we have
\[
\|\mathscr N(v)-\mathscr N(w)\|_Y\le\Cl{const:non-linear}r(3\|\bar W\|_X^2+3r\|\bar W\|_X+r^2)\|v-w\|_X,
\]
where
\[
\Cr{const:non-linear}=4\|\bar W\|_\infty(c_\infty^3+2c_\infty^2c_4^2).
\]
\end{proposition}

\begin{proof}
We have
\[
\big\|(a_{\sigma(1)}\cdot a_{\sigma(2)})(a_{\sigma(3)}\times\Delta a_{\sigma(4)})\big\|_{L^2}\le\|a_{\sigma(1)}\cdot a_{\sigma(2)}\|_{L^\infty}\|a_{\sigma(3)}\|_{L^\infty}\|\Delta a_{\sigma(4)}\|_{L^2}\le c_\infty^3\prod_{j=1}^4\|a_j\|_X,
\]
where in the first step we pulled out the maxima;
in the second we used Cauchy--Schwarz and \cref{prop:sobolev}. Analogously one finds that the sum over $i$ in the definition of $\mathcal T$ is bounded in $L^2$ by $c_\infty^2c_4^2\prod_{j=1}^4\|a_j\|_X$. 
Hence,
\[
\|\mathcal T(a_1,a_2,a_3,a_4)\|_Y\le\frac{\Cr{const:non-linear}}4\prod_{j=1}^4\|a_j\|_X.
\]
Telescoping as
\[
\mathcal T(v,v,v,v)-\mathcal T(w,w,w,w)=\mathcal T(v-w,v,v,v)+\mathcal T(w,v-w,v,v)+\mathcal T(w,w,v-w,v)+\mathcal T(w,w,w,v-w),
\]
applying the preceding bound, and the analogous two- and three-term telescopes to the other terms of $\mathscr N$, give
\[
\|\mathscr N(v)-\mathscr N(w)\|_Y\le\frac{\Cr{const:non-linear}}4(12r\|\bar W\|_X^2+12r^2\|\bar W\|_X+4r^3)\|v-w\|_X.
\]
Factoring out $4r$ proves the claim.
\end{proof}

\subsection{Energy estimate}

\begin{proposition}\label{prop:energy}
Let $W,\bar W\in X$ with $\|W-\bar W\|_X\le r$, and suppose that $u=W/|W|$ is harmonic. Assume
\[
\min_{S^3}|\bar W|-c_\infty r
\ge\Cr{const:energyparameter}>c_\infty r.
\]
Then
\begin{align*}
E\left(\frac{\bar W}{|\bar W|}\right)-\frac{c_5^2r^2}{2\Cr{const:energyparameter}^2}\left(1+\frac{4c_\infty\|\bar W\|_X}{\Cr{const:energyparameter}}\right)^2&\le E(u),\\
E(u)&\le\frac{E\left(\frac{\bar W}{|\bar W|}\right)(\Cr{const:energyparameter}+c_\infty r)^2}{\Cr{const:energyparameter}^2-3c_\infty^2r^2+2c_\infty r\Cr{const:energyparameter}}.
\end{align*}
\end{proposition}

\begin{proof}
We now bound the energy of $u$ by that of $\bar u=\bar W/|\bar W|$. Write $v=W-\bar W$ and $w=\bar u-u$. Then, since $u$ is harmonic,
$$E(\bar u)=E(u)+\tfrac12\|dw\|^2_{L^2} -\tfrac12\int|du|^2|w|^2.$$
Suppose that  $\|w\|_{C^0}^2 < 1$. Rearranging the above slightly, we can obtain two inequalities
$$E(\bar u)-\tfrac12\|dw\|_{L^2}^2\;\le\;E(u)\;\le\;\frac{E(\bar u)}{1-\|w\|_{C^0}^2}.$$

To prove the proposition, there are now two steps:
\begin{enumerate}
\item Obtain an upper bound on $\|w\|_{C^0}^2$. To do this, we use the elementary inequality $|a/|a|-b/|b||\le2|a-b|/|a|$, to get that
\[
\|w\|_{C^0}
\le\frac{2c_\infty r}{\min_{S^3}|\bar W|}
\le\frac{2c_\infty r}
{\Cr{const:energyparameter}+c_\infty r}<1.
\]
This validates the assumption that $\|w\|_{C^0}^2<1$.    \item The lower bound. Let $\varphi(X) = X/|X|$ denote the radial map on $X \in \mathbb R^3\setminus\{0\}$. In the following, we will use the following elementary estimates for $\varphi$: 
    \begin{align} \label{eq:radialbound}
    |d\varphi_X|_{\mathrm{op}} \leq \frac{1}{|X|},& &|d^2\varphi_X|_{\mathrm{op}} \leq \frac{4}{|X|^2}.
    \end{align}
    Now, to obtain a pointwise bound on $dw$, we write $$dw=d\varphi_{\bar W}(d\bar W)-d\varphi_W(d\bar W+dv)=-\,d\varphi_W(dv)+\big[d\varphi_{\bar W}-d\varphi_W\big](d\bar W) .$$ 
    To bound the first term, then by \eqref{eq:radialbound} and the inequality  $\Cr{const:energyparameter}\leq |\bar W| - |v| \leq |W|$, we have
    \begin{align} \label{eqn:ubound}
|d\varphi_W(dv)|\le\frac{|dv|}{|W|}\le\frac{|dv|}{\Cr{const:energyparameter}}.    
\end{align}
    For the second term, by the fundamental theorem of calculus, we have that
    $$d\varphi_W-d\varphi_{\bar W}=\int_0^1 d^2\varphi_{\bar W+tv}(v,\cdot)\,dt.$$
    Then, as $\Cr{const:energyparameter} \leq |\bar W| - |v| \leq |\bar W+tv |$, for all $t \in [0,1]$, using the second estimate in \eqref{eq:radialbound}, we have
    \begin{align}
     |\big[d\varphi_{\bar W}-d\varphi_W\big](d\bar W)| \leq \frac{4\,|v|\,|d\bar W|}{\Cr{const:energyparameter}^2}.   
    \end{align}
    Putting these both together, and integrating over $S^3$, we get 
    $$\|dw\|_{L^2}\le\frac{\|dv\|_{L^2}}{\Cr{const:energyparameter}}+\frac{4 \| v \|_{C^0}\|d\bar W\|_{L^2}}{\Cr{const:energyparameter}^2}\le\frac{c_5 r}{\Cr{const:energyparameter}} \left( 1+\frac{4\,c_\infty \,\|\bar W\|_{X}}{\Cr{const:energyparameter}} \right),$$ where we have used Proposition \ref{prop:sobolevderivatives} in the second inequality.
    \qedhere
\end{enumerate}

\end{proof}

\section{Conclusion of the proof of the main theorems}
\label{section:conclusion}

\begin{proof}[Proof of \cref{thm:main}]
\emph{Existence and Hopf invariant.}
For every $q =2,3,4,5,6$, the estimates in \cref{prop:injectivity-estimate,prop:approximate-bounds,prop:nonlinear} verify the hypotheses of \cref{thm:NK} with $M$, $\varepsilon = \Cr{const:residual}$, $r$ from \cref{master-table} and
\[
k=\Cr{const:non-linear}r(3\|\bar W\|_X^2+3r\|\bar W\|_X+r^2).
\]
Hence, there exists $v\in X$ with $\|v\|_X\le r$ and $F(v)=0$.

Set $W=\bar W+v$ and let $V = V_q$ be the reference map from \eqref{eq:referencemap}. Since $|V|\le1$ by \cref{lem:homotopysup}, we have
\begin{equation}\label{eq:homotopyineq}
W\cdot V\ge\bar W\cdot V-|v||V|\ge\Cr{const:homotopy}-c_\infty r>0,
\end{equation}
and $\min|\bar W|\ge\Cr{const:homotopy}>c_\infty r$. Thus \cref{prop:zeros} gives a smooth $G$-equivariant harmonic map $u=W/|W|$. Since $u\cdot V>0$, \cref{prop:homotopy} gives $H(u)=H(V/|V|) = q$.

\emph{Energy.}
The certified bound in \cref{master-table} gives $\Cr{const:energyparameter}>c_\infty r$. Applying \cref{prop:energy} and the approximate energy bounds in \cref{prop:approximate-bounds} gives $\Cr{const:energy-lower2}\le E(u)\le\Cr{const:energy-upper2}$, where
\begin{align*}
\Cl{const:energy-lower2}&:=\Cr{const:energy-lower}-\frac{c_5^2r^2}{2\Cr{const:energyparameter}^2}\left(1+\frac{4c_\infty\|\bar W\|_X}{\Cr{const:energyparameter}}\right)^2,\\
\Cl{const:energy-upper2}&:=\frac{\Cr{const:energy-upper}(\Cr{const:energyparameter}+c_\infty r)^2}{\Cr{const:energyparameter}^2-3c_\infty^2r^2+2c_\infty r\Cr{const:energyparameter}}.
\end{align*}
The resulting numerical bounds are given in \cref{master-table}. These bounds exclude $\tfrac{E(u)}{8\pi^2}$ from being an integer, so $u$ cannot have the form $u = R\circ\pi\circ A$ by \cref{prop:hopfenergy}.
\end{proof}

\begin{proof}[Proof of \cref{thm:index}]
Given the nine test fields $\xi_a \in H^1(S^3, \R^3)$ from \cref{prop:indexbound}, we let $\delta$ be as in \cref{prop:quadraticestimate}. Evaluating the formula for $\delta$ in interval arithmetic, using $r$ and $\Cr{const:energyparameter}$ from \cref{master-table}, the bound $\|\bar W\|_X\le\Cr{const:Wbar-X}$ from \cref{prop:approximate-bounds}, and the exact constants $c_\infty,c_4,\beta$ of \cref{prop:sobolev}, gives $\delta\le2.798\cdot10^{-4}$, see also \cref{master-table} so $\lambda^*+\delta\le-0.04109+0.00028<0$.
Since the fields are orthogonal for $Q_{\bar W}$ and the $H^1$-inner product, the inequality $\lambda^*+\delta<0$ makes $G+\delta H$ negative definite. Thus \cref{lem:negdef,prop:index-equivalence} give $\operatorname{Ind}Q_u\ge9$.
\end{proof}

\section*{Acknowledgements}
\noindent EH is supported by the S\~ao Paulo Research Foundation (FAPESP).
This work was supported in part by a grant of access to OpenAI models through the ChatGPT for Academic Researchers program.

\noindent DP was supported by an Imperial College Research Fellowship, partly funded by Schmidt Futures.

\noindent HSE is supported by FAPESP grants \mbox{[2021/04065-6]}, associated with \emph{BRIDGES: Brazil--France interplays in Gauge Theory, extremal structures and stability}; and \mbox{[2024/00923-6]}, associated with \emph{CBG: Brazilian Centre for Geometry}. HSE is also supported by the Brazilian National Council for Scientific and Technological Development (CNPq), through Research Productivity Fellowship \mbox{[307145/2025-5]}, level PQ-A.

\noindent JS is supported by FAPESP grant \mbox{[2023/02809-3]}, linked to the BRIDGES collaboration \mbox{[2021/04065-6]}. Their work in this project was not supported by OpenAI. 

\section*{Data availability}
Self-contained code for the numerics used in this proof, as well as certificates for proof steps, are available at this paper's respective GitHub repository: \href{https://github.com/edhirst/eells_counterexamples}{\texttt{github.com/edhirst/eells\_counterexamples}}

\section*{AI usage}
This work started with an application of physics-informed neural networks to model the map between $S^3$ and $S^2$; aligning with ideas introduced in \cite{Hirst:2026qwi}.
In this construction, to ensure the Hopf invariant was fixed throughout the learning process, the neural architecture modelled a homotopy map from the standard Hopf-related construction in each class at $t=0$ to an open endpoint at $t=1$.
This AI architecture was then trained with the goal of minimising the harmonic energy at the $t=1$ endpoint to identify potentially new harmonic maps. 
Training produced many promising candidates across Hopf invariant classes, and studying these approximations led to the identification of the symmetry used in the numerical approximant function class used in this proof.

With this information, Claude (Anthropic) then assisted with building the full candidate function, writing code for the numerics, and drafting key stages of this computer-assisted proof. Codex (OpenAI) was later used to assess and review what was produced. This was reviewed again, and then written up entirely by the authors, who retain full responsibility for the mathematical content and computational results.

\appendix

\section{Numerical constants}
In this appendix, we list the numerical constants needed for the proof of the main theorems in \S \ref{section:conclusion}. We include the choice of cut-off $l$ of \S\ref{subsection:linear-estimate} into high and low modes for the Laplacian. As required by \S \ref{section:conclusion} for \cref{thm:main}, in every column $Mk<1$, $M(\varepsilon+kr)\le r$, $\Cr{const:homotopy}-c_\infty r>0$ and $\|w\|_{C^0}<1$, moreover, for the index computation   $\lambda^*+\delta<0$, see \cref{prop:quadraticestimate}. The constants $c_\infty,c_4,c_5,\beta$ are the exact values coming from our Sobolev norms in  \cref{prop:sobolev,prop:sobolevderivatives}; the certificates use interval enclosures of them, which lie in the intervals shown.
\begin{table}[H]
\centering\small
\begin{tabular}{lccccc}
\toprule
Hopf invariant $q$ & $2$ & $3$ & $4$ & $5$ & $6$\\
\midrule
$\deg\bar W$ & $40$ & $80$ & $112$ & $152$ & $184$\\
$l$ & $40$ & $40$ & $60$ & $60$ & $80$\\
$\Cr{const:residual}$ & $\le 1.987\cdot 10^{-8}$ & $\le 5.180\cdot 10^{-11}$ & $\le 3.855\cdot 10^{-11}$ & $\le 4.657\cdot 10^{-11}$ & $\le 6.482\cdot 10^{-11}$\\
$\Cr{const:energy-lower}$ & $\ge 125.2744885849$ & $\ge 155.6036384307$ & $\ge 177.7689807885$ & $\ge 195.1505049029$ & $\ge 209.4236634609$\\
$\Cr{const:energy-upper}$ & $\le 125.2744885857$ & $\le 155.6036384337$ & $\le 177.7689807937$ & $\le 195.1505049116$ & $\le 209.4236634726$\\
$\Cr{const:homotopy}$ & $\ge 0.2499$ & $\ge 0.2450$ & $\ge 0.1917$ & $\ge 0.1347$ & $\ge 0.08728$\\
$\Cr{const:non-linear}$ & $\le 0.03534$ & $\le 0.03534$ & $\le 0.03534$ & $\le 0.03534$ & $\le 0.03534$\\
$\Cr{const:energy-lower2}$ & $\ge 125.2744885848$ & $\ge 155.6036384307$ & $\ge 177.7689807885$ & $\ge 195.1505049029$ & $\ge 209.4236634609$\\
$\Cr{const:energy-upper2}$ & $\le 125.2744885858$ & $\le 155.6036384337$ & $\le 177.7689807937$ & $\le 195.1505049116$ & $\le 209.4236634726$\\
$\Cr{const:energyparameter}$
& $0.99999986$ & $0.99999999$ & $0.99999999$
& $0.99999999$ & $0.99999998$\\
$\Cr{const:Wbar-X}$
& $69.058$ & $99.765$ & $130.98$ & $162.42$ & $193.99$\\
$\Cr{const:Wbar-C0}$
& $1.00001$ & $1.00001$ & $1.00001$ & $1.00001$ & $1.00001$\\
$\eta$
& $2.338\cdot10^{-9}$ & $1.161\cdot10^{-8}$
& $4.118\cdot10^{-8}$ & $1.207\cdot10^{-7}$
& $2.097\cdot10^{-7}$\\
$b_1$ & $\le 11.41$ & $\le 17.17$ & $\le 22.93$ & $\le 28.69$ & $\le 34.46$\\
$b_0$ & $\le 5.751$ & $\le 9.562$ & $\le 16.45$ & $\le 25.50$ & $\le 36.64$\\
$\gamma$ & $\ge 1278$ & $\ge 1032$ & $\ge 2405$ & $\ge 2039$ & $\ge 3861$\\
$\|A^{-1}\|$ & $\le 0.1476$ & $\le 0.1363$ & $\le 0.1276$ & $\le 0.1225$ & $\le 0.1192$\\
$\|B\|$ & $\le 160.7$ & $\le 261.0$ & $\le 490.3$ & $\le 648.1$ & $\le 994.1$\\
$\|C\|$ & $\le 138.6$ & $\le 181.3$ & $\le 373.8$ & $\le 433.8$ & $\le 715.5$\\
$\|A^{-1}B\|$ & $\le 0.1038$ & $\le 0.1754$ & $\le 0.1511$ & $\le 0.2067$ & $\le 0.1793$\\
$\|\mathscr L^{-1}\|_{L^2\to L^2}$ & $\le 0.1521$ & $\le 0.1482$ & $\le 0.1353$ & $\le 0.1380$ & $\le 0.1298$\\
$\|\mathscr L^{-1}\|_{Y\to X}$ & $\le M=23.67$ & $\le M=48.63$ & $\le M=77.67$ & $\le M=122.7$ & $\le M=165.6$\\
$r$ & $9.5\cdot 10^{-7}$ & $5.1\cdot 10^{-9}$ & $6\cdot 10^{-9}$ & $1.2\cdot 10^{-8}$ & $2.2\cdot 10^{-8}$\\
$Mk$ & $\le 0.0114$ & $\le 2.62\cdot 10^{-4}$ & $\le 8.48\cdot 10^{-4}$ & $\le 0.00412$ & $\le 0.0146$\\
$M(\varepsilon+kr)$ & $\le 4.82\cdot 10^{-7}$ & $\le 2.52\cdot 10^{-9}$ & $\le 3\cdot 10^{-9}$ & $\le 5.76\cdot 10^{-9}$ & $\le 1.11\cdot 10^{-8}$\\
$\|w\|_{C^0}$ & $\le 2.69\cdot 10^{-7}$ & $\le 1.45\cdot 10^{-9}$ & $\le 1.7\cdot 10^{-9}$ & $\le 3.4\cdot 10^{-9}$ & $\le 6.23\cdot 10^{-9}$\\
$\delta$ & $\le 2.798\cdot 10^{-4}$ & -- & -- & -- & --\\
$\lambda^*$ & $\le -0.04109$ & -- & -- & -- & --\\
\midrule
$c_\infty$ & \multicolumn{5}{c}{$\sqrt{\frac{\pi^2/6-5/4}{2\pi^2}}\in[0.1414482009672520,\,0.1414482009672527]$}\\[2pt]
$c_4$ & \multicolumn{5}{c}{$\left(\frac{2}{9\pi^2}\right)^{1/4}\in[0.3873663890421022,\,0.3873663890421032]$}\\[2pt]
$c_5$ & \multicolumn{5}{c}{$\frac{2\sqrt2}{9}\in[0.3142696805273541,\,0.3142696805273548]$}\\[2pt]
$\beta$ & \multicolumn{5}{c}{$(2\pi^2)^{-1/4}\in[0.4744249983287935,\,0.4744249983287952]$}\\[2pt]
\bottomrule
\end{tabular}
\caption{Constants of the certificate for the approximate solutions $\bar W=\bar W_q$ of Hopf invariant $q=2,\dots,6$ in  \cref{thm:main}. Upper bounds are rounded up and lower bounds down. }
\label{master-table}
\end{table}

\section{The approximate solutions}

The approximate solutions $\bar W=\bar W_q$ are given by equivariant polynomial maps $S^3 \rightarrow \R^3$ in the span of the basis of \cref{prop:harmonic-bases}, of the degree given in \cref{master-table}. Their leading coefficients are listed in \cref{tab:coefficients}.

In the table below, we distinguish the two types of fields in the basis of \cref{prop:harmonic-bases},  which we will call type $U$, $V$: 
\begin{align*}
(\Re Y_{q(n-2),n,j},\Im Y_{q(n-2),n,j},0)& &\sqrt2(0,0,\Re Y_{qn,n,j})    
\end{align*}
with $n$ even, and $n\geq1$ odd, respectively. 

For every $q$, more than $99.99\%$ of the squared $L^2$-norm of the approximate solution $\|\bar W_q\|_{L^2}^2$ comes from the fields of type $U$ with $n=0,2$ and of type $V$ with $n=1$. 

\begin{table}[H]
\centering\small
\begin{tabular}{crrrrrrr}
\toprule
type & $n$ & $j$ & $q=2$ & $q=3$ & $q=4$ & $q=5$ & $q=6$\\
\midrule
U & $0$ & $0$ & $-1.525663$ & $-1.023630$ & $-0.744217$ & $-0.571476$ & $-0.456312$\\
U & $0$ & $1$ & $\phantom{-}0.291036$ & $\phantom{-}0.343169$ & $\phantom{-}0.327313$ & $\phantom{-}0.295575$ & $\phantom{-}0.262980$\\
U & $0$ & $2$ & $-0.063119$ & $-0.122853$ & $-0.149767$ & $-0.156980$ & $-0.154470$\\
U & $0$ & $3$ & $\phantom{-}0.015490$ & $\phantom{-}0.046593$ & $\phantom{-}0.070918$ & $\phantom{-}0.085308$ & $\phantom{-}0.092249$\\
U & $0$ & $4$ & $-0.004064$ & $-0.018564$ & $-0.034639$ & $-0.047342$ & $-0.055932$\\
U & $0$ & $5$ & $\phantom{-}0.001097$ & $\phantom{-}0.007679$ & $\phantom{-}0.017388$ & $\phantom{-}0.026785$ & $\phantom{-}0.034396$\\
U & $0$ & $6$ & $-0.000301$ & $-0.003259$ & $-0.008930$ & $-0.015421$ & $-0.021433$\\
U & $0$ & $7$ & $\phantom{-}0.000083$ & $\phantom{-}0.001407$ & $\phantom{-}0.004670$ & $\phantom{-}0.009015$ & $\phantom{-}0.013518$\\
\addlinespace[2pt]
U & $2$ & $0$ & $\phantom{-}3.478069$ & $\phantom{-}3.719205$ & $\phantom{-}3.794535$ & $\phantom{-}3.822473$ & $\phantom{-}3.834396$\\
U & $2$ & $1$ & $-0.898244$ & $-1.325303$ & $-1.500491$ & $-1.576199$ & $-1.611709$\\
U & $2$ & $2$ & $\phantom{-}0.205195$ & $\phantom{-}0.529395$ & $\phantom{-}0.730052$ & $\phantom{-}0.837857$ & $\phantom{-}0.895588$\\
U & $2$ & $3$ & $-0.049602$ & $-0.209652$ & $-0.367253$ & $-0.476414$ & $-0.544815$\\
U & $2$ & $4$ & $\phantom{-}0.012938$ & $\phantom{-}0.083949$ & $\phantom{-}0.185306$ & $\phantom{-}0.275383$ & $\phantom{-}0.341920$\\
U & $2$ & $5$ & $-0.003515$ & $-0.034444$ & $-0.093954$ & $-0.159794$ & $-0.216736$\\
U & $2$ & $6$ & $\phantom{-}0.000971$ & $\phantom{-}0.014507$ & $\phantom{-}0.048120$ & $\phantom{-}0.093011$ & $\phantom{-}0.137848$\\
U & $2$ & $7$ & $-0.000271$ & $-0.006241$ & $-0.024974$ & $-0.054421$ & $-0.087877$\\
U & $2$ & $8$ & $\phantom{-}0.000076$ & $\phantom{-}0.002725$ & $\phantom{-}0.013140$ & $\phantom{-}0.032072$ & $\phantom{-}0.056196$\\
U & $2$ & $9$ & $-0.000021$ & $-0.001201$ & $-0.007000$ & $-0.019059$ & $-0.036091$\\
U & $2$ & $10$ & $\phantom{-}0.000006$ & $\phantom{-}0.000533$ & $\phantom{-}0.003766$ & $\phantom{-}0.011423$ & $\phantom{-}0.023300$\\
U & $2$ & $11$ & $-0.000002$ & $-0.000237$ & $-0.002042$ & $-0.006901$ & $-0.015129$\\
\addlinespace[2pt]
U & $4$ & $0$ & $-0.012727$ & $-0.008865$ & $-0.005399$ & $-0.003338$ & $-0.002147$\\
\addlinespace[2pt]
V & $1$ & $0$ & $\phantom{-}2.021640$ & $\phantom{-}1.464655$ & $\phantom{-}1.081160$ & $\phantom{-}0.823037$ & $\phantom{-}0.644811$\\
V & $1$ & $1$ & $-0.519577$ & $-0.637786$ & $-0.612132$ & $-0.547986$ & $-0.479881$\\
V & $1$ & $2$ & $\phantom{-}0.124740$ & $\phantom{-}0.261339$ & $\phantom{-}0.321739$ & $\phantom{-}0.334523$ & $\phantom{-}0.324242$\\
V & $1$ & $3$ & $-0.032001$ & $-0.106537$ & $-0.166037$ & $-0.199211$ & $-0.212729$\\
V & $1$ & $4$ & $\phantom{-}0.008680$ & $\phantom{-}0.044216$ & $\phantom{-}0.085636$ & $\phantom{-}0.117712$ & $\phantom{-}0.137966$\\
V & $1$ & $5$ & $-0.002410$ & $-0.018788$ & $-0.044528$ & $-0.069560$ & $-0.089123$\\
V & $1$ & $6$ & $\phantom{-}0.000675$ & $\phantom{-}0.008145$ & $\phantom{-}0.023429$ & $\phantom{-}0.041290$ & $\phantom{-}0.057582$\\
V & $1$ & $7$ & $-0.000190$ & $-0.003580$ & $-0.012478$ & $-0.024677$ & $-0.037306$\\
V & $1$ & $8$ & $\phantom{-}0.000054$ & $\phantom{-}0.001587$ & $\phantom{-}0.006717$ & $\phantom{-}0.014861$ & $\phantom{-}0.024273$\\
V & $1$ & $9$ & $-0.000015$ & $-0.000708$ & $-0.003646$ & $-0.009017$ & $-0.015874$\\
V & $1$ & $10$ & $\phantom{-}0.000004$ & $\phantom{-}0.000317$ & $\phantom{-}0.001992$ & $\phantom{-}0.005507$ & $\phantom{-}0.010436$\\
\addlinespace[2pt]
V & $3$ & $0$ & $-0.011423$ & $-0.007876$ & $-0.004977$ & $-0.003238$ & $-0.002198$\\
\midrule
\multicolumn{3}{l}{$\deg\bar W_q$} & $40$ & $80$ & $112$ & $152$ & $184$\\
\multicolumn{3}{l}{number of coefficients} & $217$ & $630$ & $969$ & $1482$ & $1847$\\
\bottomrule
\end{tabular}
\caption{Leading coefficients $a_k$ of the approximate solutions $\bar W_q=\sum_k a_k\Psi_k$ in the orthonormal basis of Proposition \ref{prop:harmonic-bases}. Listed are all basis fields whose coefficient has modulus at least $0.01$ for some $q$, rounded to 6 decimals. The complete lists, with the exact IEEE double coefficients used by the certificate, are the files \texttt{data/Wbar\_q$q$\_D$D$.json} of the repository.}
\label{tab:coefficients}
\end{table}

\section{The test fields for the index}
\label{appendix:index-fields}
In this appendix, we give a list of nine negative directions $\xi_a$ of the quadratic form $Q_{\bar{W}}$ in \eqref{eq:quadraticform} for our approximate solution $\bar W$ with Hopf invariant $2$, which were found using computer-assisted methods. 

To understand these nine directions somewhat systematically, we use the induced action \eqref{eqn:symmetries} of the subgroup $G':=U(1) \times \mathbb Z_2 \subset G_{q=2}$ generated by $e^{i\theta}$ and $e^{i\pi}\varepsilon$ on $ H^1(S^3,\mathbb R^3)$. By the Peter-Weyl theorem, this induces a splitting (over $\C$) of $ H^1(S^3,\mathbb R^3)$ into isotypic $G'$-representations, orthogonally with respect to the $H^1$ inner product. Moreover, since $\bar W$ is $G'$-equivariant, $Q_{\bar W}$ is $G'$-invariant, and so distinct isotypic components are also orthogonal for $Q_{\bar W}$. 

We can write a basis of isotypic $G'$-representations explicitly over $\C$. For $Y=Y_{m,n,j}$
put
\begin{equation}\label{eq:basisPMT}
P_{m,n,j}=\tfrac1{\sqrt2}(Y,-iY,0),\qquad M_{m,n,j}=\tfrac1{\sqrt2}(Y,iY,0),\qquad T_{m,n,j}=(0,0,Y).
\end{equation}
By Lemma \ref{lem:action} the action of $G'$ on this basis is given by 
\begin{align*}
e^{i\theta} \cdot \xi=e^{i\nu\theta}\xi,& & (e^{i\pi}\varepsilon)\cdot\xi=p\,\xi,   
\end{align*}
where $\nu=2n-m-4,\ 2n-m+4,\ 2n-m$ and
$p=(-1)^l,\ (-1)^l,\ (-1)^{l+1}$ for each of $P,M,T$ respectively.

Let $S_{\nu,p}$ be the closed complex linear span of $P,M,T$ in \eqref{eq:basisPMT} with weights $(\nu,p)$ for the $G'$-action, i.e. an isotypic component for the $G'$-action on $ H^1(S^3,\C^3)$. Complex conjugation maps $S_{\nu,p}$ to $S_{-\nu,p}$. \Cref{tab:index-fields} lists the five fields $f_1,\dots,f_5$ for $q=2$, which give the nine negative directions $\xi_1=f_1$, $\Re f_a$, $\Im f_a$ ($a=2,\dots,5$).

For $f\in S_{\nu,p}$ with $\nu\ne0$, rotation by $\pi/(2\nu)$
sends $(\Re f,\Im f)$ to $(-\Im f,\Re f)$.
Hence, $\Re f$ and $\Im f$ are orthogonal and have equal quadratic values for $Q_{\bar{W}}$.
Hence, the nine given fields are nonzero and pairwise orthogonal.

\begin{table}[H]
\centering\footnotesize\setlength{\tabcolsep}{3pt}
\begin{tabular}{lr@{\hspace{12pt}}lr@{\hspace{12pt}}lr@{\hspace{12pt}}lr@{\hspace{12pt}}lr}
\toprule
\multicolumn{2}{c}{$f_{1}\in S_{0,-}$} & \multicolumn{2}{c}{$f_{2}\in S_{1,-}$} & \multicolumn{2}{c}{$f_{3}\in S_{2,-}$} & \multicolumn{2}{c}{$f_{4}\in S_{2,+}$} & \multicolumn{2}{c}{$f_{5}\in S_{3,-}$}\\
\midrule
$T_{0,0,0}$ & $-100$ & $T_{1,1,0}$ & $-89$ & $M_{0,-1,0}$ & $100$ & $T_{0,1,0}$ & $100$ & $M_{1,0,0}$ & $100$\\
$T_{0,0,1}$ & $-56$ & $M_{-1,-2,0}$ & $17$ & $T_{-2,0,0}$ & $41$ & $M_{2,0,0}$ & $-66$ & $T_{-1,1,0}$ & $-76$\\
$M_{2,-1,0}$ & $39$ & $M_{3,0,0}$ & $100$ & $M_{0,-1,1}$ & $-48$ & $T_{0,1,1}$ & $1$ & $M_{1,0,1}$ & $-15$\\
$P_{-2,1,0}$ & $39$ & $P_{-1,2,0}$ & $16$ & $P_{0,3,0}$ & $-64$ & $M_{-2,-2,0}$ & $-22$ & $T_{-1,1,1}$ & $-11$\\
$T_{0,0,2}$ & $12$ & $T_{-3,-1,0}$ & $46$ & $T_{-2,0,1}$ & $-35$ & $M_{2,0,1}$ & $-2$ & $M_{-3,-2,0}$ & $23$\\
$M_{-2,-3,0}$ & $10$ & $T_{1,1,1}$ & $5$ & $T_{2,2,0}$ & $-42$ & $P_{-2,2,0}$ & $-23$ & $M_{1,0,2}$ & $-9$\\
$M_{2,-1,1}$ & $-2$ & $M_{-1,-2,1}$ & $-14$ & $M_{0,-1,2}$ & $19$ & $T_{-4,-1,0}$ & $-25$ & $P_{-3,2,0}$ & $33$\\
$P_{-2,1,1}$ & $-2$ & $M_{3,0,1}$ & $-8$ & $M_{4,1,0}$ & $27$ & $T_{0,1,2}$ & $-8$ & $P_{1,4,0}$ & $-18$\\
$P_{2,3,0}$ & $10$ & $P_{-5,0,0}$ & $-63$ & $P_{-4,1,0}$ & $29$ & $M_{-2,-2,1}$ & $8$ & $T_{-5,-1,0}$ & $30$\\
$T_{-4,-2,0}$ & $9$ & $P_{-1,2,1}$ & $-14$ & $P_{0,3,1}$ & $27$ & $M_{2,0,2}$ & $3$ & $T_{-1,1,2}$ & $12$\\
$T_{0,0,3}$ & $1$ & $T_{-3,-1,1}$ & $-17$ & $T_{-2,0,2}$ & $14$ & $P_{-6,0,0}$ & $26$ & $T_{3,3,0}$ & $-8$\\
$T_{4,2,0}$ & $9$ & $T_{1,1,2}$ & $3$ & $T_{2,2,1}$ & $16$ & $P_{-2,2,1}$ & $8$ & $M_{-3,-2,1}$ & $-8$\\
$M_{-2,-3,1}$ & $-4$ & $M_{-1,-2,2}$ & $5$ & $M_{0,-1,3}$ & $-6$ & $P_{2,4,0}$ & $-1$ & $M_{1,0,3}$ & $5$\\
$M_{2,-1,2}$ & $-1$ & $P_{-5,0,1}$ & $14$ & $M_{4,1,1}$ & $-8$ & $T_{-4,-1,1}$ & $7$ & $M_{5,2,0}$ & $3$\\
$M_{6,1,0}$ & $-6$ & $P_{-1,2,2}$ & $5$ & $P_{-4,1,1}$ & $-8$ & $T_{0,1,3}$ & $3$ & $P_{-7,0,0}$ & $-31$\\
$P_{-6,-1,0}$ & $-6$ & $P_{3,4,0}$ & $1$ & $P_{0,3,2}$ & $-8$ & $T_{4,3,0}$ & $-1$ & $P_{-3,2,1}$ & $-11$\\
$P_{-2,1,2}$ & $-1$ & $T_{-3,-1,2}$ & $4$ & $T_{-2,0,3}$ & $-4$ & $M_{-2,-2,2}$ & $-2$ & $P_{1,4,1}$ & $6$\\
$P_{2,3,1}$ & $-4$ & $T_{1,1,3}$ & $-1$ & $T_{2,2,2}$ & $-4$ & $M_{2,0,3}$ & $-1$ & $T_{-5,-1,1}$ & $-7$\\
$T_{-4,-2,1}$ & $-2$ & $T_{5,3,0}$ & $1$ &  &  & $M_{6,2,0}$ & $1$ & $T_{-1,1,3}$ & $-4$\\
$T_{0,0,4}$ & $-1$ &  &  &  &  & $P_{-6,0,1}$ & $-5$ & $T_{3,3,1}$ & $2$\\
$T_{4,2,1}$ & $-2$ &  &  &  &  & $P_{-2,2,2}$ & $-2$ &  & \\
\midrule
\multicolumn{2}{c}{$\le -0.4693$} & \multicolumn{2}{c}{$\le -0.04850$} & \multicolumn{2}{c}{$\le -0.05998$} & \multicolumn{2}{c}{$\le -0.2957$} & \multicolumn{2}{c}{$\le -0.04109$}\\
\bottomrule
\end{tabular}
\caption{The test fields of \cref{prop:indexbound}: $f_a=\sum c\,\Phi$ with Gaussian-integer coefficients $c$ in the basis $\Phi\in\{P_{m,n,j},M_{m,n,j},T_{m,n,j}\}$ of \eqref{eq:basisPMT} in the isotypic component $S_{\nu,p}$. They are rounded eigenvectors of the matrix given by $Q_{\bar W}$ restricted to modes of degree $\le 8$ of each isotypic component. The last row is the certified bound on $Q_{\bar W}(\xi)/\|\xi\|_{H^1}^2$ for $\xi=f_1$, respectively $\xi=\Re f_a,\Im f_a$ ($a\ge2$); its maximum is $\lambda^*$ listed in table \ref{master-table}.}
\label{tab:index-fields}
\end{table}

\section{Exact quadrature}

Throughout the paper, we computed many integrals of spherical harmonics. In this section, we use the techniques of exact quadrature to carry out computation of these integrals. For the convenience of the reader, we give a brief summary of this standard technique, see also \cite{DavisRabinowitz}, \cite[Sections~1 and~3]{TrefethenExactness}.

The classical method of \emph{quadrature} replaces an integral by a finite weighted sum of values of the integrand at prescribed points, called quadrature points or nodes. A quadrature rule is \emph{exact} on a given space of functions if this finite sum equals the integral for every function in that space: the two linear functionals agree there.
For polynomials of bounded degree, suitable nodes and weights give such an identity, rather than merely an approximation as the number of nodes tends to infinity.

In our case, we are integrating functions $S^3 \rightarrow \C$ that are invariant under a circle action. This reduces the problem to two-dimensional quadrature. We record here our exact quadrature rule.
\paragraph{}
Fix integers $N,M\ge2$. For $0\le i\le N$ and $0\le j<M$, set
\[
  s_i=\frac{1-\cos(\pi i/N)}2,\qquad
  x_{ij}=\bigl(\sqrt{s_i},
    \sqrt{1-s_i}\,e^{2\pi\mathrm{i}j/M}\bigr)\in S^3.
\]
A univariate polynomial $P(s)$ of degree at most $N$ is uniquely determined by its values at the $N+1$ distinct $s_i$, since the difference of two such polynomials cannot have $N+1$ roots unless it vanishes.
Given the numerical values $s_i$, the \emph{Lagrange interpolation polynomials} $\ell_i\in\R[s]$ form a basis of polynomials of degree at most $N$ adapted to these values: these $\ell_i$ are the uniquely defined polynomials, such that each has degree at most $N$, equals $1$ at $s_i$, and $0$ at the other $s_k$-values, i.e. $\ell_i(s_k)=\delta_{ik}$. Every polynomial $P$ of degree at most $N$ satisfies $P(s)=\sum_{i=0}^N P(s_i)\ell_i(s)$; integrating this identity will give us our quadrature weights \eqref{eqn:integration-weights}.

So, to each node $x_{ij} \in S^3$ we assign the weight
\begin{align}
\label{eqn:integration-weights}
  \omega_i=\frac{2\pi^2}{M}\int_0^1\ell_i(s)\,ds,
\end{align}
independent of $j$.
For a function $f:[0,1]\times[0,2\pi]\to\C$ one can now approximate
the value of the integral
\[
\int_0^1\int_0^{2\pi}f(s,\phi)\,d\phi\,ds
\approx \sum_{i=0}^{N}\sum_{j=0}^{M-1}
\frac{\omega_i}{\pi}f(s_i,2\pi j/M).
\]
Here is what motivates the weights of this quadrature formula.
The \emph{periodic trapezoidal rule} approximates an integral over $[0,2\pi]$ by $2\pi$ times the average of the function's values at the $M$ equally spaced nodes $2\pi j/M$.
Thus every node has weight $2\pi/M$
\cite[p.~394, equation~(3.2)]{TrefethenWeideman2014}.
Then, the \emph{Clenshaw--Curtis rule} on $[0,1]$ integrates the polynomial interpolating the function at the nodes $s_i$, giving weights $\int_0^1\ell_i(s)\,ds$
\cite[p.~135]{TrefethenExactness}.

The weights $\omega_i/\pi$, with $\omega_i$ as in
\cref{eqn:integration-weights}, are precisely the product
of these two sets of weights: apply one rule in each variable and multiply their weights.
The extra factor of $\pi$ in \cref{eqn:integration-weights} adapts the quadrature for integrals on the $3$-sphere as explained in the following proposition:

\begin{proposition}\label{prop:orbit-quadrature}
Let $q\ge2$ be an integer and let $f:\R^4 \to\C$ be a polynomial of degree at most $p$, invariant under
$(z,w)\mapsto(e^{\mathrm{i}\theta}z,e^{-q\mathrm{i}\theta}w)$.
If $p\le 2N$ and $p<(q+1)M$, then
\[
  \int_{S^3}f\,d\mathrm{vol}
  =\sum_{i=0}^{N}\sum_{j=0}^{M-1}\omega_i f(x_{ij}).
\]
\end{proposition}

\begin{proof}
    In the coordinates $(z,w)=(\sqrt{s} e^{i\alpha}, \sqrt{1-s} e^{i \beta})$, the volume element becomes $\frac{1}{2} ds \, d\alpha \, d\beta$.
    Circle invariance and Fubini's theorem therefore give
    \[
    \int_{S^3}f\,d\mathrm{vol}
    =
    \pi\int_0^1\int_0^{2\pi}
    f\bigl(\sqrt{s},\sqrt{1-s}\,e^{\mathrm i\phi}\bigr)
    \,d\phi\,ds,
    \]
    where we used the substitution $\phi=q \alpha+\beta$.

    By linearity, it suffices to consider monomials of the form $f=z^a \bar z^b w^c \bar w^d$ for $a-b=q(c-d)$.
    We have
    \[
    \begin{aligned}
    f\bigl(\sqrt{s},\sqrt{1-s}\,e^{\mathrm i\phi}\bigr)
    &=(\sqrt{s})^{a+b}
      \bigl(\sqrt{1-s}\,e^{\mathrm i\phi}\bigr)^c
      \bigl(\sqrt{1-s}\,e^{-\mathrm i\phi}\bigr)^d
      =s^{(a+b)/2}(1-s)^{(c+d)/2}
      e^{\mathrm i(c-d)\phi}.
    \end{aligned}
    \]
    Here, its angular frequency $c-d$ satisfies
    \[
    |c-d|
    =
    \frac{|a-b|+|c-d|}{q+1}
    \le
    \frac{a+b+c+d}{q+1}
    \le
    \frac{p}{q+1}
    <
    M,
    \]
    where in the first step we used $a-b=q(c-d)$;
    in the second step we used the triangle inequality;
    in the third we used the assumption that $f$ has at most degree $p$;
    and in the fourth the hypothesis on $p$. The $M$-point trapezoidal rule is therefore exact in the angular variable by \cite[Corollary 3.3]{TrefethenWeideman2014}.

    The angular integral and the angular quadrature sum vanish unless $c=d$.
    In that case, we have $a=b$, and the integrand becomes $s^a(1-s)^c$ for $a+c \leq p/2 \leq N$.
    Thus, Clenshaw--Curtis quadrature at the nodes $s_i$ integrates this polynomial exactly, see \cite[Eqn.~(19.4) and Theorem~19.1]{TrefethenATAP} (after rescaling $[-1,1]$ to $[0,1]$).
    The weight in \cref{eqn:integration-weights} is $\pi$ times the product of the weights for the trapezoidal rule \cite[Eqn.~(3.2)]{TrefethenWeideman2014} and Clenshaw--Curtis quadrature from \cite[Eqn.~(19.4)]{TrefethenATAP}.
\end{proof}

\section{Visualisations of the harmonic maps}
\label{app:visualisations}

We illustrate the harmonic maps $u_q$ of \cref{thm:main} with $H(u_q)=q$, together with $q=1$, for which the same construction returns exactly the Hopf fibration $\pi$ of \cref{def:hopffibration}.
The figures are computed from $\bar W_q/|\bar W_q|$ with the scripts in \texttt{analysis/} of the repository; since $u_q=W/|W|$ with $W=\bar W_q+v$ and $\|v\|_{C^0}\le c_\infty r$ (\cref{master-table}), this differs from $u_q$ by less than $10^{-6}$ in $C^0$.
Points of $S^3$ are drawn in $\R^3$ by the stereographic projection $(z,w)\mapsto(\Re z,\Im z,\Re w)/(1-\Im w)$.

\begin{figure}[htbp]
\centering
\begin{subfigure}[t]{0.49\textwidth}
  \centering
  \includegraphics[width=\linewidth]{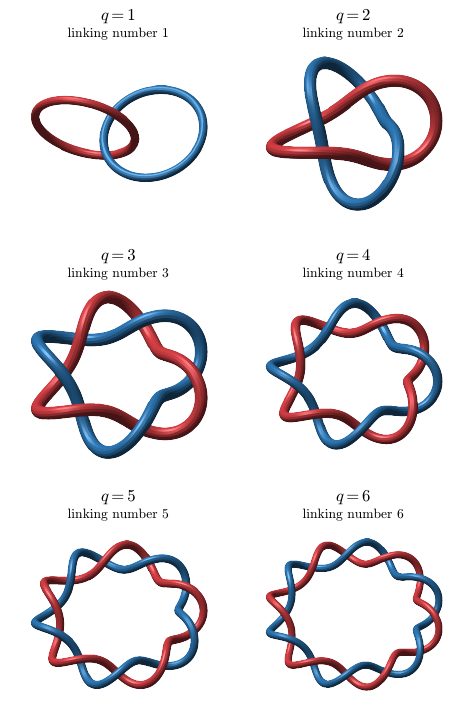}
  \caption{The pole fibres $u_q^{-1}(0,1)$ (red) and $u_q^{-1}(0,-1)$ (blue).}
  \label{fig:fibres-poles}
\end{subfigure}\hfill
\begin{subfigure}[t]{0.49\textwidth}
  \centering
  \includegraphics[width=\linewidth]{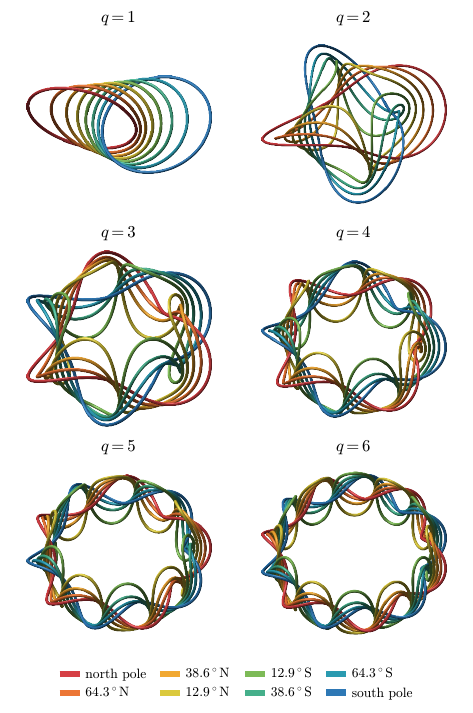}
  \caption{Fibres over eight equally spaced points of a meridian of $S^2$.}
  \label{fig:fibres-meridian}
\end{subfigure}
\caption{Fibres of $u_q$, $q=1,\dots,6$, stereographically projected to $\R^3$.}
\label{fig:fibres}
\end{figure}

\paragraph{Fibres and linking.}
In \cref{fig:fibres-poles}, the fibres of $u_q$ over the north and south poles of $S^2$ are two closed curves, and the image shows how many times one winds around the other.
That is, the absolute value of the linking number can be seen as the number of times the blue curve crosses over the red curve.
There is one crossing for the Hopf fibration ($q=1$, the Hopf link) and $q$ for $u_q$. This linking number is the classical geometric picture of the Hopf invariant \cite[\S17]{BottTu}, so the figure lets one see $H(u_q)=q$ directly.
\Cref{fig:fibres-meridian} adds the fibres over six further points of a meridian of $S^2$.
They sweep out an annulus bounded by the two pole fibres, providing a visualisation of this transition which twists $q$ times, and with any two of the fibres linking $q$ times.

\begin{figure}[htbp]
\centering
\begin{subfigure}[t]{0.49\textwidth}
  \centering
  \includegraphics[width=\linewidth]{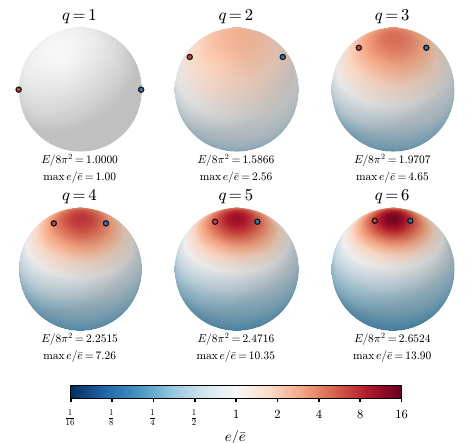}
  \caption{Energy density $e$ on the orbit sphere $S^3/S^1$, relative to its mean $\bar e=E/2\pi^2$; dots: the pole fibres.}
  \label{fig:energy-orbits}
\end{subfigure}\hfill
\begin{subfigure}[t]{0.49\textwidth}
  \centering
  \includegraphics[width=\linewidth]{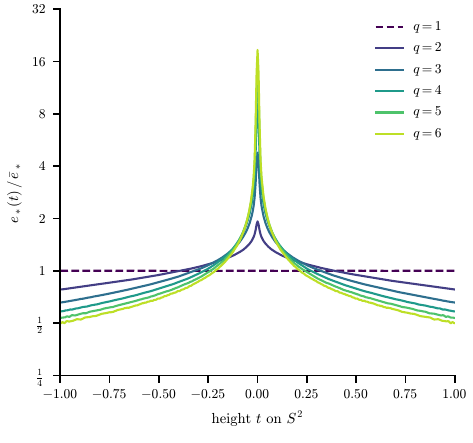}
  \caption{Energy density $e_*$ on the target $S^2$, relative to its mean $\bar e_*=E/4\pi$, against the height $t$.}
  \label{fig:energy-target}
\end{subfigure}
\caption{Energy of $u_q$, $q=1,\dots,6$.}
\label{fig:energy}
\end{figure}

\paragraph{Energy.}
Our maps $u_q$ are told apart from the maps of \cref{conj:eells} by their energy, so it is natural to ask where this energy sits.
The energy density $e=\frac12|du_q|^2$ is constant along the circle orbits of $G_q$, so \cref{fig:energy-orbits} draws it on the sphere of orbits for each $q$, with height $|z|^2-|w|^2$ and longitude $\phi = \arg(z^qw)$.
On this $2$-sphere for all $q$ equal areas correspond to equal volumes of $S^3$ since $d\mathrm{vol}=\pi\,ds\,d\phi$ with $s=|z|^2$ (see the proof of \cref{prop:orbit-quadrature}); the dots are the pole fibres of \cref{fig:fibres-poles}.
\Cref{fig:energy-target} instead shows how the energy is spread over the target $S^2$, which by symmetry depends only on the height $t$.
For the Hopf fibration both are constant; as $q$ grows, the energy concentrates near the orbit $w=0$ in $S^3$ and over the equator of $S^2$.

\newpage

\bibliographystyle{plain}
\bibliography{references}

\end{document}